\documentclass[12pt,notitlepage,twoside,a4paper]{amsart}
 \usepackage{amsfonts}

\usepackage{amsmath,amssymb,enumerate}

\usepackage{epsfig,fancyhdr,color}

\usepackage{amssymb}
\usepackage{amsmath,amsthm}
\usepackage{latexsym}
\usepackage{amscd}
\usepackage{psfrag}
\usepackage{graphicx}
\usepackage[latin1]{inputenc}
\usepackage[all]{xy}
\usepackage[mathcal]{eucal}

\definecolor{NoteColor}{rgb}{1,0,0}

\renewcommand{\textsc}{\textcolor{red}}

\newtheorem{theorem}{\rm\bf Theorem}[section]
\newtheorem{proposition}[theorem]{\rm\bf Proposition}
\newtheorem{lemma}[theorem]{\rm\bf Lemma}
\newtheorem{corollary}[theorem]{\rm\bf Corollary}
\newtheorem*{theorem 1}{\rm\bf Proposition 1}
\newtheorem*{theorem 2}{\rm\bf Proposition 2}

\theoremstyle{definition}
\newtheorem{definition}[theorem]{\rm\bf Definition}

\theoremstyle{remark}

\def\interieur#1{\mathord{\mathop{\kern 0pt #1}\limits^\circ}}

\title[On Grothendieck's construction of 
Teichm\"uller space]{On Grothendieck's construction of 
\\
Teichm\"uller space}

\makeindex 

\author[N. A'Campo, L. Ji and A. Papadopoulos]{Norbert A'Campo, Lizhen Ji \\
and Athanase Papadopoulos}

\address{N. A'Campo: Universit\"at Basel,  Mathematisches Institut, 
\\
Spiegelgasse 1, 4051 Basel, Switzerland
\\
and 
Erwin Schr\"odinger International Institute of Mathematical Physics, 
\\
Boltzmanngasse 9, 1090, Wien, Austria
\\
email:\,\tt{norbert.acampo@gmail.com}}
\address{L. Ji: Department of Mathematics, University of Michigan\\ Ann Arbor, MI 48109, USA
\\
and 
Erwin Schr\"odinger International Institute of Mathematical Physics, 
\\
Boltzmanngasse 9, 1090, Wien, Austria
\\
email:\,\tt{lji@umich.edu}}
\address{A. Papadopoulos: Institut de Recherche Math\'ematique Avanc\'ee,
\\
Universit{\'e} de Strasbourg and CNRS,\\
7 rue Ren\'e Descartes, 67084 Strasbourg Cedex, France
\\
  and Erwin Schr\"odinger International Inititute of Mathematical Physics, 
  \\
  Boltzmanngasse 9, 1090, Wien, Austria 
\\
email:\,\tt{papadop@math.unistra.fr}}

\begin{document}

\maketitle

 \begin{abstract} In his 1944 paper \emph{Ver\"anderliche Riemannsche Fl\"achen}, Teichm\"uller defined a structure of complex manifold on the set of isomorphism classes of marked closed Riemann surfaces of genus $g$. The complex manifold he obtained is the space called today  Teichm\"uller
 space. In the same paper, Teichm\"uller introduced the so-called universal Teichm\"uller curve -- a space over Teichm\"uller space where the fiber above each point is a Riemann surface representing that point. In fact, Teichm\"uller proved the existence of the Teichm\"uller curve as a space of Riemann surfaces parametrized by an analytic space, with an existence and uniqueness theorem establishing this analytic structure. This result was later reformulated and proved by Grothendieck in a series of ten lectures he gave at Cartan's seminar in 1960-1961. In his approach, Grothendieck replaced Teichm\"uller's explicit parameters by a general construction of fiber bundles whose base is an arbitrary analytic space. This work on Teichm\"uller space led him to recast the bases of analytic geometry using the language of categories and functors.  In Grothendieck's words, the Teichm\"uller curve becomes a space representing a functor from the category of analytic spaces into the category of sets. In this survey, we comment on Grothendieck's series of lectures.
 
 The survey is primarily addressed to low-dimensional topologists and geometers. In presenting Grothendieck's results, we tried to explain or rephrase in more simple terms some notions that are usually expressed in the language of algebraic geometry. However, it is not possible to short-circuit the language of categories and functors.   
 
 The survey is also addressed to those algebraic geometers who wish to know how the notion of moduli space evolved in
connection with Teichm\"uller theory.  

Explaining the origins of mathematical ideas contributes in dispensing justice to their authors  and it usually renders the theory that is surveyed more attractive. 

\end{abstract}



\medskip

\noindent {\it The final version of this paper will appear as a chapter in Volume VI of the {\it Handbook of Teichm\"uller theory}. This volume is dedicated to the memory of Alexander Grothendieck.}

\tableofcontents

\section{Introduction}

One of the last works of Teichm\"uller is his paper \cite{T32} published in 1944, in which he constructs Teichm\"uller space as a universal object and defines the Teichm\"uller (or universal) curve, and where he equips these spaces with complex-analytic structures.\index{Teichm\"uller space!complex structure} Compared to its importance, this paper remains very poorly known. A translation and a commentary were recently published, see \cite{T32} and \cite{T32C}. Teichm\"uller's remarkable result was very influential on Alexander Grothendieck, who reformulated it at the end of the 1950s and who presented in the academic year 1960-1961 a new construction of this space, in a series of ten lectures he gave at  Cartan's\footnote{Although this will seem obvious to everybody, we point out that all the occurrences of the name Cartan in this chapter refer to Henri Cartan.\index{Cartan, Henri} (This will exempt us from adding the first initial H. each time we write his name, as is usually required in the mathematical literature, to distinguish him from his father Elie Cartan.)} seminar in Paris. 
 Cartan used to ask the speakers at his seminar to provide written versions for their talks. These written texts were usually distributed as mimeographed notes at the subsequent seminar meetings.\footnote{Talking about Cartan and his seminar, Douady writes in \cite{Douady-Cartan}: ``He did not tolerate the slightest inaccuracy, the slightest imprecision, and he criticized the speaker to the point of destabilizing him. [...] But what really mattered for Cartan were the notes of the expos\'es. Again, he would tolerate no imprecision. It was out of the question to say that two groups were isomorphic without specifying an isomorphism between them, or to say that a diagram commuted up to sign: this sign had to be given. Above all he wanted the text to be perfectly clear. For this he asked you to revise the text as many times as necessary. \emph{Vingt fois sur le m\'etier remettez votre ouvrage (Redo your work twenty times)} said Boileau. With Cartan it was rather thirty than twenty times. His point was that a text should have at least thirty readers, otherwise there's no point in writing it. Therefore, if you spend half an hour to spare the reader a minute of perplexity, then it's well worth it. He would return your manuscript covered with annotations made with a red pen in his small, curly handwriting. Then you would revise it and give him back again.} Thus, we have at our disposal a series of papers by Grothendieck on Teichm\"uller theory (\cite{G1} to \cite{G10}). 
We shall present them briefly, saying a few words on each paper, at some places providing some details on definitions and important concepts.

Grothendieck formulated Teichm\"uller's theorem using the languages of categories, of algebraic geometry and of analytic geometry. It is always good to remember that Riemann surfaces themselves are at the same time analytic and algebraic objects, and therefore it is natural that they have been studied using methods of analytic as well as algebraic geometry. 
We recall in this respect that any closed Riemann surface (that is, a complex one-dimensional manifold) of genus $g\geq 2$ admits holomorphic embedding in the projective space $\mathbb{CP}^{5g-5}$, using the Riemann-Roch theorem for the third tensor power of the canonical bundle. 
By a theorem of Chow \cite{Chow}, this image  is defined by homogeneous polynomial equations. By post-composing with a generic projection we obtain an embedding in $\mathbb{CP}^{3}$. The case $g=0$ is part of Riemann's uniformization theorem. The case $g=1$ is the one of elliptic curves and uses the Weierstrass $p$- function $\wp$ to yield an embedding in 
$\mathbb{CP}^{2}$.

 In the course of his work on the subject, Grothendieck considered that classical algebraic geometry did not have enough tools for formulating and proving the existence of Teichm\"uller space as a universal object carrying a complex structure, and for giving an algebraic model for moduli spaces.  In the introduction to the first lecture \cite{G1}, he  writes:
``In doing this, the necessity of reshaping the foundations of analytic geometry, inspired by the theory of schemes, will be manifest."\footnote{[Chemin faisant, la n\'ecessit\'e deviendra manifeste de revoir les fondements de la g\'eom\'etrie analytique.]} In particular, the notion of schemes which he had newly introduced turned out to be useful in dealing with the problems of moduli of Riemann surfaces and in other moduli problems. At the same time, Grothendieck also brought new ideas into Teichm\"uller theory.

 The general title of Grothendieck's 10 papers is  \emph{Techniques de construction en g\'eom\'etrie analytique} (Construction techniques in analytic geometry), and each paper carries its own title. 

 In these papers, Grothendieck constructs Teichm\"uller space as a complex space representing a functor. Grothendieck worked essentially in the analytic setting, because the required passage to quotient by a group acting properly was available in analytic, but not in algebraic geometry.  In this analytic geometry setting, it is not unusual that the existence of some moduli space amounts to the representability of some functor. Grothendieck's papers also contain the definition of and the existence theorem for the universal Teichm\"uller curve equipped with its complex structure.  In Section \ref{s:functor} below, we shall describe at length the content of the first paper,  titled \emph{Description axiomatique de l'espace de Teichm\"uller} (Axiomatic description of Teichm\"uller space) \cite{G1}, because it contains the motivating ideas and an outline of the whole theory. We shall also describe in detail the last paper (Expos\'e X), titled \emph{Construction de l'espace de Teichm\"uller} (Construction of Teichm\"uller space), which contains the proofs of the main results.
 
  Grothendieck, in constructing Riemann's moduli space and other moduli spaces, used extensively Teichm\"uller's idea of enhancing the structures by markings in order to overcome the 
difficulties that are caused by nontrivial automorphism groups of the objects representing points in  moduli spaces. In a letter to Serre written on  November 5, 1959 (see \cite{GS-correspondence}), Grothendieck writes:
\begin{quote} \small 
I have already come to the practical conclusion that every time that my criteria show 
that no modular variety (or rather, moduli scheme) for the classification
of (global or infinitesimal) variations of certain structures (complete
non-singular varieties, vector bundles, etc.) cannot exist, despite good hypotheses of 
flatness, properness, and if necessary non-singularity, the only reason is the existence
of automorphisms of the structure which prevent the descent from
working. [...] The remedy in moduli theory seems to me to be to eliminate bothersome
automorphisms by introducing additional structures on the objects being studied: points or differential forms, etc. on the varying varieties
(a process is already used for curves), trivializations at sufficiently
many points of the vector bundles one wants to vary, etc.\footnote{[La conclusion pratique \`a laquelle je suis arriv\'e d\`es maintenant, c'est que chaque fois que, en vertu de mes crit\`eres, une vari\'et\'e de modules, (ou plut\^ot, un sh\'ema de modules) pour la classification des vari\'et\'es (globales, ou infinit\'esimales) de certaines structures (vari\'et\'es compl\`etes non singuli\`eres, fibr\'es vectoriels etc.) ne peut exister, malgr\'e de bonnes hypoth\`eses de platitude, propret\'e et non singularit\'e \'eventuellement, la raison en est seulement l'existence d'automorphismes de la structure qui emp\^eche la descente de marcher. [...] La panac\'ee en th\'eorie des modules me semble alors, si les automorphismes nous emb\^etent, de les \'eliminer en introduisant des donn\'ees suppl\'ementaires sur les objets qu'on \'etudie : points ou formes diff\'erentielles etc. sur les vari\'et\'es qu'on fait varier (proc\'ed\'e d\'ej\`a utilis\'e pour les courbes), trivialisations en suffisamment de points des fibr\'es vectoriels qu'on veut faire varier, etc.]}
\end{quote} One important consequence of this work which results in viewing the enhanced Riemann moduli space as a space representing a functor, besides providing an analytic structure on this space, is that  the same methods apply to other moduli whic can be realized as quotients of Hermitian symmetric spaces by some discrete groups.

 In this survey, while we present Grothendieck's ideas, we try to explain a few notions of algebraic geometry, relating them when possible to notions which are familiar to Teichm\"uller theorists, but without entering into the technical details, often at the expense of being too vague for an expert in algebraic geometry who would read this text. All the displayed results in the following sections (Theorem, Proposition, etc.) are, with very few exceptions, contained as such in Grothendieck's Cartan seminars.

We end this introduction by stressing the fact that the point of view of Grothendieck on defining the Teichm\"uller space as a space representing a functor not only led to the standard modern formulation
and construction of other moduli spaces in algebraic geometry and to other generalizations, but it also deepened our understanding of Teichm\"uller space and clarified some aspects of the nature of this space which is more than a space of points, but a space of Riemann surfaces and their deformations. 

\medskip

\noindent{\it Acknowledgements.---} The authors acknowledge support from the U.S. National Science Foundation grants DMS 1107452, 1107263, 1107367 ``RNMS: GEometric structures And Representation varieties" (the GEAR Network). Part of this work was done during a stay of the three authors at Galatasaray University (Istanbul) and at the Erwin Schr\"odinger Institute (Vienna). We are grateful to Pierre Deligne, Antoine Ducros, H\'el\`ene Esnault and Hiroaki Nakamura for a careful reading of a previous version of this chapter and for suggesting several improvements.

\section{A note on Alexander Grothendieck and on Cartan's seminar}

We start with a personal touch, from the first author of the present chapter.  During my undergraduate studies at the University of  Montpellier I, I
worked, in the autumn of 1962, as vintage helper, during the grape harvesting period, in the village of Vendargues.         
During a meal with the farmers -- I remember their names: the farmer Jean-Henri Teissier, his wife Annette, their son Jean-Louis  and their daughter Yvette --  and other vintage workers, the farmer, who knew that I was studying mathematics, asked me if I was familiar with a mathematician called Grothendieck. I had never heard of this name and my  answer was a clear ``no."
Two years later, my teacher, Professor Jean-Pierre Laffon, suggested  that I read the paper \emph{Sur quelques points d'alg\`ebre homologique}, by Alexander
Grothendieck, which was published in 1957 in the T\^ohoku Mathematical Journal. I was amazed to hear for the second time in France this Dutch sounding name under so
different circumstances. When I went back to the family Teissier in Vendargues, I got the following explanation: Grothendieck had lived with his mother in Mairargues, in a 
house at walking distance from the house of the farmers' family.\footnote{In \emph{R\'ecoltes et semailles}\index{Grothendieck!R\'ecoltes et semailles}\index{R\'ecoltes et semailles (Grothendieck)} \S\,2.1, Note 1, Grothendieck mentions this place (The translation from this writing of Gorthendieck's is ours): ``Between 1945 and 1948, I used to live, with my mother, in a small hamlet, at about a dozen kilometers from Montpellier, Mairargues (by Vendargues), lost among the vines. (My father had disappeared in Auschwitz, in 1942.) We were living meanly, on my student's scholarship. To help make both ends meet, I used to work in harvesting, every year, and after the vintage, in that of the ``vin de grapillage" [a non-commercial table wine] (which, as I was told, was illegal...) Furthermore, there was a garden which, without any effort from our part, provided us with figs, spinach, and (towards the end) tomatoes,  planted by an indulgent neighbor, right in the middle of an ocean of magnificent poppies. This was the life! But sometimes too short in the joints, when it came to replacing a spectacle frame, or a pair of threadbare shoes. Fortunately, for my mother, weak and ill after her long stay in the concentration camps, we were entitled to free medical assistance. We would never had been able to pay any doctor..." [Entre 1945 et 1948, je vivais avec ma m\`ere dans un petit hameau \`a une dizaine de kilom\`etres de Montpellier, Mairargues (par Vendargues), perdu au milieu des vignes. (Mon p\`ere avait disparu \`a Auschwitz, en 1942.) On vivait chichement sur ma maigre bourse d'\'etudiant. Pour arriver \`a joindre les deux bouts, je faisais les vendanges chaque ann\'ee, et apr\`es les vendanges, du vin de grapillage, que j'arrivais \`a \'ecouler tant bien que mal (en contravention, para\^\i t-il, de la l\'egislation en vigueur. . . ) De plus il y avait un jardin qui, sans avoir \`a le travailler jamais, nous fournissait en abondance figues, \'epinards et m\^eme (vers la fin) des tomates, plant\'ees par un voisin complaisant au beau milieu d'une mer de splendides pavots. C'\'etait la belle vie - mais parfois juste aux entournures, quand il s'agissait de remplacer une monture de lunettes, ou une paire de souliers us\'es jusqu'\`a la corde. Heureusement que pour ma m\`ere, affaiblie et malade \`a la suite de son long s\'ejour dans les camps, on avait droit \`a l'assistance m\'edicale gratuite. Jamais on ne serait arriv\'es \`a payer un m\'edecin. . .]}

Cartan's seminar lasted 16 years, from 1948 to 1964. The principle was the following: every year, a  new theme was chosen, and the seminar talks were supposed to start from scratch and give complete proofs. The seminar was attended by several talented mathematicians and it certainly was a major element in the fact that Paris became at that time the world mathematical capital. Another major element was IH\'ES, the Institut des Hautes \'Etudes Scientifiques, founded in 1958. Cartier writes in \cite{Cartier-Gro} that ``Grothendieck dominated the first ten years of the institute," and Serre writes in \cite{Serre-Cartan} that the ``impressive" \emph{S\'eminaire de g\'eom\'etrie alg\'ebrique} of Grothendieck at IH\'ES, which lasted 10 years (1960-1969), was the successor of Cartan's seminar.\index{Cartan seminar} The ten seminar talks that Grothendieck gave at the Cartan seminar in 1960-1961 were a very important step in his major work on the foundations of algebraic geometry.  Hubbard recalls, in \cite{Hubbard-A}, that Douady, who was present at this series of talks, told him: ``Adrien described to me Cartan's resentment when Grothendieck started his first talk saying `In this talk and in the next ones, I will present ...'; Cartan felt that his seminar did not belong to him any more."\footnote{[Adrien m'a d\'ecrit le d\'epit de Cartan lorsque Grothendieck a commenc\'e son premier expos\'e disant :  Dans cet expos\'e et les suivants j'exposerai ... , et que Cartan a senti son s\'eminaire lui \'echapper.]}

   In  his personal and mathematical autobiography\emph{R\'ecoltes et semailles},\index{Grothendieck!R\'ecoltes et semailles}\index{R\'ecoltes et semailles (Grothendieck)} (1986) Grothendieck, at several times, mentions Cartan's seminar, as well as his own seminar, which he started in the 1960s, and which was known under the nickname SGA (S\'eminaire de G\'eom\'etrie Alg\'ebrique),\footnote{We learned from Deligne that the name ``S\'eminaire de G\'eom\'etrie Alg\'ebrique du Bois Marie" came later, and that what is now called SGA3 was called SGAD (D for Demazure) and SGA4: SGAA (A for Artin).} reporting on the differences between the two seminars, but also the common points. In \S\,18, he writes: ``That which was common to the two seminars seems to me more important, and above all, that which was their essential function, their raison d'\^etre. As a matter of fact, I see two points. One function of these seminars, which is close to Bourbaki's purpose, was to prepare and to provide to everybody texts which are easily accessible (I mean, essentially complete), developing in an elaborate manner important themes which are difficult to access.\footnote{Grothendieck's footnote: ``Difficult to access, either because these themes remained imperfectly understood, or because they were known only to rare initiates, and such that the scattered relevant publications gave an inadequate image of them.
["D'acc\`es difficile," soit parce que ces th\`emes restaient imparfaitement compris, soit qu'ils n'\'etaient connus que de rares initi\'es, et que les publications \'eparpill\'ees qui en traitaient n'en donnaient qu'une image inad\'equate.]} The other function of these two seminars is to provide a place where young and motivated researchers were sure, even without being  geniuses, to be able to learn the craft of a mathematician for what concerns present date questions, in contact with eminent and caring men. To learn the craft -- that is, to put the shoulder to the wheel, and henceforth, to find the occasion of being known."\footnote{[Plus important me para\^\i t ce qui \'etait commun aux deux s\'eminaires, et surtout, ce qui me semble avoir \'et\'e leur fonction essentielle, leur raison d'\^etre. \`A vrai dire j'en vois deux. Une des fonctions de ces s\'eminaires, proche du propos de Bourbaki, \'etait de pr\'eparer et de mettre \`a la disposition de tous des textes ais\'ement accessibles (j'entends, essentiellement complets), d\'eveloppant de fa\c con circonstanci\'ee des th\`emes importants et d'acc\`es difficile. L'autre fonction de ces s\'eminaires, \'etait de constituer un lieu o\`u des jeunes chercheurs motiv\'es \'etaient s\^urs, m\^eme sans \^etre des g\'enies, de pouvoir apprendre le m\'etier de math\'ematicien sur des questions de pleine actualit\'e, au contact d'hommes \'eminents et bienveillants. Apprendre le m\'etier -- c'est \`a dire, mettre la main \`a la p\^ate, et par l\`a-m\^eme, trouver l'occasion de se faire conna\^\i tre.]}

One may also add that Grothendieck, before these ten lectures at Cartan's seminar, had already given lectures on moduli. Mumford, in his paper \cite{MG}, remembers Grothendieck's lectures at Harvard, in 1958, two years before the Cartan seminar lectures. He writes: ``My involvement came about because I had been studying the construction of varieties classifying families of algebraic structures, especially moduli spaces of vector bundles and of curves. Whereas I had thought loosely of such a classifying space as having a `natural' one-one correspondence with the set of objects in question (just as Riemann and Picard had), Grothendieck expressed it with functors. This was clearly the right perspective. There were `fine' moduli spaces\index{fine moduli space}\index{moduli space!fine} which carried a universal family of objects, e.g. a universal family of curves from which all other families were unique pull-backs. Therefore they represented the functor of all such families. And there were also `coarse' moduli spaces,\index{coarse moduli space}\index{moduli space!coarse} the best possible representable approximation to the desired functor (the approximation being caused e.g. by the fact that some curves had automorphisms)."

Several comments on the life and work of Alexander Grothendieck are made in other chapters of this volume, \cite{Po} and \cite{ACJP}.

\section{An introduction to the major ideas}\label{s:ideas}
 
The results to which Grothendieck refers in his work on Teichm\"uller space and which may have been his motivation in giving these lectures at Cartan seminar are contained in one of the the last papers that Teichm\"uller published \cite{T32}. Let us start by recalling the content of that paper. This is the paper in which Teichm\"uller equipped the space that carries his name\footnote{The name ``Teichm\"uller space" was given by Andr\'e Weil.} with a complex structure and computed its dimension, thus giving a precise meaning to Riemann's moduli problem and solving it. This complex-analytic structure is characterized by a certain universal property, and Teichm\"uller constructed a fiber bundle  over that space, which was called later on the \emph{Teichm\"uller universal curve}.  This construction is different from the construction of Teichm\"uller space as a metric space, which was given by Teichm\"uller in his 1939 paper \cite{T20} (see also the commentary \cite{T20C}). 

Teichm\"uller's paper \cite{T32} contains several new ideas, including the following:
\begin{enumerate}
\item The idea of rigidifying Riemann surfaces (or non-singular algebraic curves) by introducing markings, for the construction of  a non-singular moduli space.
\item The idea of a \emph{fine}\index{fine moduli space}\index{moduli space!fine} moduli space. This is a moduli space together with a universal family over it, such that every family of curves is a pull-back from that universal family by a uniquely determined map from the base into the moduli space. 
\item The idea that Teichm\"uller space is defined by the functor it represents. It is contained in the definition given by Teichm\"uller of Teichm\"uller space as a complex space satisfying a universal property. This idea was developed and generalized to other moduli spaces by Grothendieck, for whom the idea of a universal property should be replaced by that of the representability of a certain functor.
\item The existence and uniqueness of the universal Teichm\"uller curve. This curve, after Teichm\"uller and Grothendieck, was rediscovered by Ahlfors and by Bers. The existence and uniqueness result introduced at the same time the first fiber space over Teichm\"uller space. 

In fact, Teichm\"uller's paper is the first paper where families of Riemann surfaces are studied.
\item The proof of the fact that the automorphism group of the universal Teichm\"uller curve is the extended mapping class group.
\item The idea of using the period map to define the complex structure of Teichm\"uller space. 
\item A precise meaning to Riemann's heuristic count of the ``number of moduli". The $3g-3$ moduli that Riemann announced become, in Teichm\"uller's work, a complex dimension.
\end{enumerate}
The first five items in this list are investigated again in the papers by Grothendieck that we are reporting on here.

Teichm\"uller's  statement in \cite{T32} is the following:
\begin{theorem} \label{th:T}
There exists an essentially unique globally analytic family of topologically determined Riemann surfaces $\underline{\mathfrak{H}}[\frak{c}]$, where $\mathfrak{c}$ runs over a $\tau$-dimensional complex analytic manifold $\mathfrak{C}$ such that for any topologically determined Riemann surface $\underline{\mathfrak{H}}$ of genus $g$ there is one and only one $\mathfrak{c}$ such that the Riemann surface  $\underline{\mathfrak{H}}$ is conformally equivalent to an $\underline{\mathfrak{H}}[\frak{c}]$ and such that the family $\underline{\mathfrak{H}}[\frak{c}]$ satisfies the following universal property: If $\underline{\mathfrak{H}}[\frak{p}]$ is any globally analytic family of topologically determined Riemann surfaces with base $\mathfrak{B}$, there is a holomorphic map $f:\mathfrak{B}\to\mathfrak{C}$ such that the family $\underline{\mathfrak{H}}[\frak{p}]$ is the pull-back by $f$ of the family $\underline{\mathfrak{H}}[\frak{c}]$. 
\end{theorem}

A few remarks are in order. In this statement, $\mathfrak{C}$ is the space that we call today Teichm\"uller space and $\underline{\mathfrak{H}}[\frak{c}]$ is a fiber bundle over $\mathfrak{C}$, $\underline{\mathfrak{H}}[\frak{c}]\to  \mathfrak{C}$, where the fiber above each point in $\mathfrak{C}$ is a marked Riemann surface representing this point.  The expression ``topologically determined," in Teichm\"uller's vocabulary, means,  ``marked."  The marking\index{marking} here is a topological enhancement of the structure which distinguishes it from other surfaces which are conformally equivalent to it. In general in Teichm\"uller's work, several sorts of markings are used, and not only the notion of a marking as we intend it today (that is, a homotopy class of homeomorphisms between the given Riemann surface and a fixed topological surface).\footnote{The work ``marking"\index{marking} in this chapter does not refer to the choice of  a marked point, as this term is sometimes used in algebraic geometry.} We shall see below examples of markings.\footnote{Teichm\"uller was aware of the equivalence of the various notions of markings. In his paper \cite{T29}  (see also  the commentary \cite{T29C}), he uses three sorts of markings: homotopy classes of homeomorphisms, isotopy classes of homeomorphisms, and the choice of a basis of the fundamental group.\index{marking!basis for the fundamental group} These equivalences are deep theorems in the topology of surfaces, and Teichm\"uller attributes them to Mangler \cite{Mangler}.} Finally, $\tau$ is the complex dimension, given by  $\tau= 0$ if $g=0$, $1$ is $g=1$ and  $3(g-1)$  if $g>1$. This is Riemann's ``number of moduli." The fiber bundle  $\underline{\mathfrak{H}}[\frak{c}]\to  \mathfrak{C}$ is the Teichm\"uller universal curve.

 Thus, the theorem says that if $\underline{\mathfrak{H}}[\frak{p}]\to\mathfrak{B}$ is any analytic  family of marked Riemann surfaces, then it is obtained from the universal Teichm\"uller curve $\underline{\mathfrak{H}}[\frak{c}]\to  \mathfrak{C}$ as a pull-back by some holomorphic map $f:\mathfrak{B}\to\mathfrak{C}$.
From the context, and stated in modern terms, the essential uniqueness of  $\underline{\mathfrak{H}}[\frak{c}]\to  \mathfrak{C}$  announced at the beginning of the statement of the theorem means that the family is unique up  to the action of the mapping class group. In other words, the map is unique up to the choice of a marking.

For an English translation and a commentary of Teichm\"uller's paper we refer the reader to \cite{T32} and \cite{T32C}.

 Teichm\"uller's paper is difficult to read, both for analysts and for low-dimensional topologists and geometers, because of its very concise style, and also because it uses heavily the language of algebraic geometry (function fields, divisors, valuations, places, etc.)  Grothendieck took up the idea and developed it using a new language of algebraic geometry. In turn, Grothendieck's papers are difficult to read for people not used to his language. His point of view is based on the fact that Teichm\"uller space  represents the functor of marked families of projective 
algebraic
curves.  Let us quote Grothendieck's statement:

 \begin{theorem}[Theorem 3.1. of \cite{G1}] \label{the:GG} 
 There exists an analytic space $T$ and a $\mathcal{P}$-algebraic curve $V$ above $T$ which are universal in the following sense: For every $\mathcal{P}$-algebraic curve $X$ above an analytic space $S$, there exists a unique analytic morphism $g$ from $S$ to $T$ such that $X$ (together with its $\mathcal{P}$-structure) is isomorphic to the pull-back of $V/T$ by $g$.
\end{theorem}

The reader should notice the intertwining of the words ``algebraic" and  ``analytic" in the expression  : ``algebraic curve over an analytic space." This is a family of algebraic curves which depends analytically on a parameter. In this statement, the analytic space $T$ is Teichm\"uller space and the $\mathcal{P}$-algebraic curve $V$ above $T$ is the universal curve. The term ``$\mathcal{P}$-algebraic" refers to a rigidification\index{rigidification} of the curves, and $\mathcal{P}$ is a functor, the so-called rigidifying functor.\index{functor!rigidifying}

There is another formulation by Grothendieck of this existence theorem, which is important, and which uses the language of functors. We shall give it below, but we start by recalling some basic definitions on categories and functors that the reader might need in order to understand the statements. The amount of knowledge needed here is very small.

A category $\mathcal{C}$\index{category} is defined by :
\begin{enumerate}
\item A class whose elements are called the objects of $\mathcal{C}$.
\item For each pair $(a,b)$ of elements of $\mathcal{C}$, there is a set\footnote{Although these questions will not be of interest for us here, we note that the morphisms between two elements form a set in the classical sense.} $\mathrm{Hom}(a,b)$ whose elements are called the \emph{morphisms}, or \emph{arrows}, from $a$ to $b$. Given an element $f$ of $\mathrm{Hom}(a,b)$,  the object $a$ is called the \emph{domain} and $b$ the \emph{codomain} of $f$. The expression ``morphism $f$ from $a$ to $b$" is denoted by $f:a\to b$.
\item For each triple $(a,b,c)$ of elements of $\mathcal{C}$, there is a binary operation $\mathrm{Hom}(a,b)\times \mathrm{Hom}(b,c) \to \mathrm{Hom}(a,c)$ called \emph{composition} of arrows and denoted by $(f,g)\to gf$.

\item The morphisms and the composition operation are required to satisfy the following properties:
\begin{itemize}
\item For any distinct pairs $(a,b)$ and $(a',b')$ the sets of morphisms $\mathrm{Hom}(a,b)$ and $\mathrm{Hom}(a',b')$ are disjoint.
\item (Associativity) For any four objects $a$, $b$, $c$ and $d$ in $\mathcal{C}$ and for any three morphisms $f:a\to b$, $g:b\to c$ and $h:c\to d$, we have
\[h(gf)=(hg)f.\]
\item (Identity element) For any object $a$ in $\mathcal{C}$, there is a morphism $\mathrm{id}_a:a\to a$ satisfying $f \mathrm{id}_a=f$ for every morphism $f$ with domain $a$ and $\mathrm{id}_a g=g$ for every morphism $g$ with codomain $a$.
\end{itemize}

A morphism $f:a\to b$ is said to be invertible if there is a morphism $g:b\to a$ such that $gf=\mathrm{id}_a$ and $fg=\mathrm{id}_b$. An invertible morphism is also called an isomorphism.

\end{enumerate}
Categories existed before Gorthendieck;
they were introduced by Cartan and Eilenberg in the setting of homological algebra \cite{CE}. It was Grothendieck who used them with all their wealth.

There are standard examples of categories, namely, the category of sets and mappings, the category of groups and group homomorphisms, the category of topological spaces and continuous maps, etc. Below, we shall mention several categories that appear in Grothendieck's work. 

Given a category $\mathcal{C}$, one defines its \emph{dual category}\index{category!dual}\index{dual category}  (also called \emph{opposite category})\index{category!opposite}\index{opposite category}  $\mathcal{C}^o$ as the category which has the same objects but whose morphisms are reversed.  Thus,  $\mathcal{C}^o$ has the same class of objects and morphisms, but for any morphism of $f$ of $\mathcal{C}$, its domain in $\mathcal{C}^o$ is the codomain of $f$ in $\mathcal{C}$ and vice versa. Composition $fg$ of two morphisms in $\mathcal{C}^o$ is defined as composition $gf$ in $\mathcal{C}$.

Given two categories $\mathcal{C}$ and $\mathcal{E}$, a \emph{covariant functor}\index{covariant functor}\index{functor}\index{functor!covariant}  $F: \mathcal{C}\to \mathcal{E}$ is a function which associates to every object $a$ in $\mathcal{C}$ an object $F(a)$ in $\mathcal{E}$ and to each morphism $f:a\to b$ in $\mathcal{C}$ a morphism $F(f):F(a)\to F(b)$ in $\mathcal{E}$ satisfying the two identities, concerning composition and the identity morphism:

\begin{itemize}
\item $F(Id_a)= id_{F(a)}$ for any object $a$ in $\mathcal{C}$;

\item $F(gf)=F(g)F(f)$ for any morphisms $f:a\to b$ and $g:b\to c$ between objects $a,b,c$ in $\mathcal{C}$.
\end{itemize}

A \emph{contravariant functor}\index{contravariant functor}\index{functor!contravariant}  $F: \mathcal{C}\to \mathcal{E}$ is a covariant functor $\mathcal{C}^o\to \mathcal{E}$ such that $F(gf)=F(f)F(g)$ whenever the composition $gf$ is defined in $\mathcal{C}$.

In category theory, the notion of representable functor\index{functor!representable} is closely related to the notion of solution of a universal problem.\index{universal problem!solution}
Grothendieck gives in \cite{G4} a short introduction to representable functors.  Teichm\"uller space is one of the first interesting examples -- and may be the first one after the elementary examples (projective spaces, Grassmannians, etc.) --   of an analytic space representing a functor.\footnote{\label{f:D} Regarding the introduction of representable functors, we quote Deligne, from a letter he sent us on July 26, 2015: ``The principle that a scheme $S$ should be understood in terms of the functor $h_S:X\mapsto \mathrm{Hom}(X,S)$ owes a lot, I think, to a conversation with Yoneda. The geometry is obscured when a scheme is viewed as a ringed space, but clearer on $h_S$. For instance, if $S$ is an algebraic group, the underlying set is not a group, but $h_S$ is a functor in groups. The underlying set does not see nilpotents, while $h_S$ does. I am not sure whether the use of $h_S$ to understand $S$ comes before, or not, the work on Teichm\"uller space."} Working in the category of fiber spaces is a tool to define a functor $F$, and in this setting only isomorphisms really matter. The notion of ``fibered category" is required to define pull-backs. Since these notions are important elements in the theory, we review them in some detail.

 Fibered categories occur in the setting of fiber bundles, vector bundles or sheaves over a given topological space.\footnote{The notion of sheaf over a topological space was introduced in 1950 by Leray in \cite{Leray}. Cartan's seminar  for the year 1950-1951 was devoted to the theory of sheaves. Analytic varieties were then included in the theory of ringed topological spaces, with their underlying sheaf of holomorphic function (This is mainly due to Oka and Cartan). It was in 1954, thanks to the work of Serre \cite{Serre1954}, that this theory also included algebraic varieties. On p. 455 of \emph{R\'ecoltes et semailles}\index{Grothendieck!R\'ecoltes et semailles}\index{R\'ecoltes et semailles (Grothendieck)} (\S\,2.10), Grothendieck declares that this paper of Serre, which, as he says, is one of the few papers he read, has had an enormous influence on him; the complete quote is included in the present volume, cf. \cite{ACJP}.} For instance, in such a category, the objects are pairs $(X,E)$ where $X$ is a topological space and $E$ a fiber bundle over $X$. Morphisms are the fiber maps between the vector bundles. This category is fibered over the category of topological spaces. The functoriality amounts to the compatibility of inverse image operations with the composition of such maps. Grothendieck introduced fibered categories in his work on the descent technique \cite{Gro-Bou1} \cite{Gro-Bou2}. This can be considered as a general ``gluing" technique in topology  adapted to the setting of algebraic geometry.

A \emph{representable functor}\index{functor!representable} is a functor from a given category into the category of sets which is isomorphic to a functor of a special form, which we now define more precisely.

Let $C$ be a category and let $X$ be an object in $C$. We consider a contravariant functor $h_X$ from $C$ to the category $\mathbf{(Ens)}$ of sets\footnote{$\mathbf{(Ens)}$ stands for the French word ``Ensembles."} defined by the formula 
\[h_X(Y)=\mathrm{Hom}(Y,X)\] at the level of objects, and as follows at the level of arrows: If $Y$ and $Z$ are objects in $C$ and  $f:Y\to Z$ a morphism between them, then the image of $f$ is the map $p\mapsto p\circ f$ from the set $\mathrm{Hom}(Z,X)$ to the set $\mathrm{Hom}(Y,X)$. 
The functor $X\to h_X$ is natural with respect to the operation of taking projective limits.

\begin{definition}[Representable functor]\index{functor!representable} Let $C$ be now a category and let \[F:C\to \mathbf{(Ens)}\] be a contravariant functor from $C$ to the category of sets.
We say that $F$ is \emph{representable} if there exists an object $X$ in $C$ such that the functor $F$ is isomorphic to the functor $h_X$.
\end{definition}
We say that $F$ is \emph{represented} by the object $X$. 

(In the case where $C$ is a covariant functor, then $F$ is said to be representable if there exists an object in $C$ such that the functor $F$ is isomorphic to the functor 
$Y\mapsto \mathrm{Hom}(X,Y).$)

According to Grothendieck in \cite{Gro-Bou2},  the ``solution of a universal problem" always consists in showing that a certain functor from a certain category to the category  $\mathbf{(Ens)}$ is representable.\footnote{In the letter dated July 26, 2015 which we already mentioned (Footnote \ref{f:D}), Deligne notes the following: ``The idea that interesting objects are solutions of universal problems is older than Grothendieck and explicit in Bourbaki Ens Ch. IV (except that Bourbaki has more in mind corepresentable functors [examples: tensor product of modules, free objects of various kinds, ...] and that values of functors are viewed as sets of maps, making $E(f)$ a composition.) But Bourbaki does not tell that $S\mapsto h_S$ is fully faithful. He only tells that $h_S$ determines $S$ up to a unique isomorphism."} In principle, the fact that a functor defined over a category $C$ is representable makes this functor computable in terms of operations belonging to the category $C$.

Let us note that the general theory of representable functors applies to categories of schemes, and that when a scheme represents a functor, then it is unique as such (a result in category theory known as Yoneda's lemma).\index{Yoneda lemma} The fact that a functor is representable makes the scheme that represents it a fine moduli space for the corresponding moduli problem, that is, a moduli space carrying a universal family. Conversely, by Yoneda's lemma, a fine moduli space is unique up to isomorphism.\footnote{In fact, Yoneda's lemma is stronger than this statement. It allows the embedding of an arbitrary category into a category of functors. To an object $S$ in a category $\mathcal{C}$ is associated the contravariant functor $h_S:\mathcal{C}\to \mathbf{(Ens)}$ where $\mathbf{(Ens)}$ is the category of sets, 
which takes the value $\mathrm{Hom}_{\mathcal{C}}(X,S)$ on an arbitrary object $X$. Yoneda's lemma says that this functor $S\mapsto h_S$ is faithful.\index{Yoneda lemma}
 Thus, Yoneda's lemma replaces the study of a category by that of all functors of that category into the category of sets with functions as morphisms. In particular, the lemma tells us that in order to define a morphism between two objects $S\to S'$, it suffices to define a morphism of functors $h_S\to h_{S'}$.}

Roughly speaking, in the language of functors, the main result presented in Grothendieck's series of lectures \cite{G1} to \cite{G10} is that the functor defined on the category of analytic spaces sending a space $X$ to the set of isomorphism classes of bundle maps $E\to X$ is representable and that Teichm\"uller space is the analytic space which represents this functor.
We shall be more precise below.

 Teichm\"uller space as introduced by Teichm\"uller is probably the first non-trivial example of a complex space representing a functor. The definition of Grothendieck fits with this early example of Teichm\"uller.

The \emph{Teichm\"uller functor}\index{Teichm\"uller functor}\index{functor!Teichm\"uller} is of the form 
\[S\mapsto  \mathrm{\ isomorphism \ classes  \ of \ marked\ } X/S.\]
In \cite{G1} p. 7, Grothendieck considers the functor
\[X/S\mapsto \mathrm{\ the \ set \ of\  rigidifications}\]
and he considers the rigidifying functors of the following type:
for some group $\gamma$,
\[X/S\mapsto \mathrm{\ a \ } \gamma-\mathrm{torsor \ over  \ }X.\]
The Teichm\"uller functor becomes universal in the sense that any rigidifying functor can be deduced from the Teichm\"uller rigidifying functor (\cite{G1} p. 7) by push-out of a $\gamma$-torsor by $\gamma\to\gamma'$ and $\gamma$ in that case is what Grothendieck calls the Teichm\"uller modular group (the mapping class group).
 Now we can state Grothendieck's result  in terms of functors. We use the above notation.

 \begin{theorem}[\cite{G1} p. 8] \label{the:Functor} 

  The rigidified Teichm\"uller functor, $\mathcal{P}$,  for curves of genus $g$ is representable.
\end{theorem}

Theorems \ref{the:GG}  and \ref{the:Functor}  are equivalent.
We already recalled that according to Grothendieck, the best formulation of a universal property consists in stating that a certain functor is representable. We can quote here Grothendieck from his   Bourbaki seminar \cite{Gro-Bou2} that he gave in 1960 (the year where he gave his Cartan seminar talks): ``[...] The fact is at the basis of the \emph{notion of ``solution of a universal problem"}, such a problem always consisting in examining whether a given functor from $\mathcal{C}$ into $\mathbf{(Ens)}$ is representable."\footnote{[[...] Ce fait est \`a la base de la \emph{notion de ``solution d'un probl\`eme universel"}, un tel probl\`eme consistant toujours \`a examiner si un foncteur donn\'e  de $\mathcal{C}$ dans $\mathbf{(Ens)}$ est repr\'esentable.]}

 In the last paper in the series (Expos\'e X), the universal problem is solved. Grothendieck  constructs Riemann's moduli space. By that time, Mumford had developed the study of moduli space in the setting of geometric invariant theory, and the Riemann moduli space appears as a quasi-projective variety.

 We now make a few remarks on the comparison between the   proofs by Teichm\"uller and by Grothendieck of the existence of Teichm\"uller space as an analytic space.  Grothendieck shows, following Teichm\"uller, that by adding the marking on the fibers, which  gives rise to a covering of the space, moduli space becomes smooth. Teichm\"uller space appears as an orbifold universal covering of moduli space. (It is known that the orbifold character disappears up to a finite covering.) Teichm\"uller first constructs Teichm\"uller space, and then moduli space. He does not construct moduli space directly. Grothendieck builds moduli spaces of marked Riemann surfaces as a space representing a functor in the algebraic category, so that the space is canonically equipped with an algebraic structure (and a fortiori an analytic structure). Teichm\"uller used a topological marking of Riemann surfaces in order to construct Teichm\"uller space. Grothendieck introduced algebraic markings (a marking by cubical differentials\index{marking!cubic differentials} and a marking through level structures\index{marking!level structures}  on the first homology). Note that Teichm\"uller space is an infinite cover of moduli space, and the algebraic structure of moduli space does no lift to Teichm\"uller space. 
  
 The question of what is an algebraic deformation was one of the main questions formulated in a clear language by Grothendieck. Motivated by Teichm\"uller space, Grothendieck proposed a general approach to the construction of moduli spaces
   of algebraic varieties, in particular the Hilbert scheme. In fact, 
the year he gave his Cartan seminar series of talks on the construction of Teichm\"uller space, Grothendieck gave a series of talks at the Bourbaki seminar. In one of them, he introduced the notion of Hilbert scheme\index{Hilbert scheme}\index{scheme!Hilbert} \cite{Gro-BouIV}. This is a scheme which is a parameter space for the closed subschemes of some fixed projective space (or more generally of a projective scheme). Hilbert schemes are the building blocks of Grothendieck's theory of families of algebraic varieties. A Hilbert scheme represents the so-called Hilbert functor.\index{Hilbert functor}\index{functor!Hilbert} 
  
Let us also note that in the following year, at the Bourbaki seminar \cite{Gro-Pi}, Grothendieck sketched his theory of the Picard scheme,\index{Picard scheme}\index{scheme!Picard} based on the techniques he developed for the construction of the Hilbert scheme in  \cite{Gro-BouIV} and the techniques of passage to the quotients developed in \cite{Gro-quotients}. This is a ``representable-functor version" of the Picard group. We recall that the (absolute) Picard group of a ringed space $X$ (a topological space  equipped with a structure of sheaf of rings) is the group of isomorphism classes of invertible sheaves on that space, the group operation being tensor product. In the case where $X$ is a projective variety over an algebraically closed filed $k$, the Picard group of $X$ underlies a natural $k$-scheme called the \emph{Picard scheme}.\index{Picard scheme}\index{scheme!Picard} The paper \cite{Kleiman} by S. L. Kleiman contains an exposition of the Picard scheme, with a fascinating historical introduction.

  Teichm\"uller's idea of marking in the sense of a ``rigidification" allowed Grothendieck to remove nontrivial automorphism groups of the varieties under consideration and to construct fine moduli spaces (moduli spaces carrying universal bundles). Rigidification has the effect of a desingularization of moduli spaces by looking at a smooth covering space or a smooth variety which is mapped onto these moduli spaces, cf. \cite{G1} to \cite{G10} and \cite{Gro-Bou1},  \cite{Gro-Bou2}.

In return, Grothendieck applied ideas from algebraic geometry  to give a new point of view on the construction of Teichm\"uller space and the Teichm\"uller curve, equipped with their complex structures and he reformulated the theory in his own language.  He obtained statements that are more general than those of Teichm\"uller, and Riemann's moduli space $\mathcal{M}_{g,n}$ became the moduli stack of algebraic curves of genus $g$ with $n$ distinguished points. 

Let us be more explicit. Paraphrasing and expanding the introduction in \cite{G1}, Grothendieck's goal in his series of lectures, as announced in the first lecture, is the following:

\begin{enumerate}
\item To introduce a general functorial mechanism for a global theory of moduli. Teichm\"uller theory is one example to which this formalism applies, but the theory also applies for instance to families of elliptic curves (this is the case of genus 1), which  so far (according to Grothendieck) had not been made very explicit in the literature.\footnote{The case of genus 1 is treated in Teichm\"uller's paper \cite{T32}. Also, in the paper \cite{T20}, Teichm\"uller showed that Teichm\"uller space, with the Teichm\"uller metric, is isometric to the hyperbolic plane.}
\item To give a ``good formulation" of a certain number of moduli problems for analytic spaces. We recall that the expression ``moduli problem" has its origin in Riemann's observation that the conformal type of a closed Riemann surface of genus $g\geq 2$ depends on $3g-3$ complex parameters, and that the goal of the so-called ``theory of moduli" is to make such a statement precise and to describe such parameters. In fact, Grothendieck gave a precise formulation of several moduli problems besides the one of Riemann surfaces, e.g. moduli of Hilbert schemes of points or Hilbert schemes of subvarieties, or Hilbert embeddings with a given Hilbert polynomial.   these moduli problems were used as step-stones, but they also have an independent interest.  According to Grothendieck, the state of the art in most of the situations is such that one could only ``conjecture some reasonable existence theorems."

\item Under some ``projectivity  hypothesis"\footnote{This projectivity condition is satisfied for holomorphic line bundles with enough sections, since this allows to separate points and tangent vectors. In the non-singular case, it is equivalent to admitting a K\"ahler structure whose K\"ahler class is integral (in $H^2(X,\mathbb{Z})$). The condition is necessary for the question of representability.} for the morphisms that will be considered, to give some existence theorems for the problems in (2). \emph{``This will include in particular the existence theorem for Teichm\"uller space."}\footnote{Grothendieck underlines.}
\item \label{44} Grothendieck says that for this purpose, it will be necessary to reconsider the foundations of analytic geometry by getting inspiration from the theory of schemes. In particular, it will be important to admit nilpotent elements\index{nilpotent element} in the local rings defining analytic spaces, and also in more general spaces that consist in families of spaces,  so that the theorems are stated with all their strength.

\end{enumerate} 

To understand Item (\ref{44}), recall that a classical\footnote{In analytic and algebraic geometry, the term ``classical" often refers to the pre-Grothendieck transformation of these fields.} complex analytic space is defined as a ``ringed space,"\index{ringed space} that is, a system of local rings (or a \emph{sheaf}\index{sheaf}) of holomorphic functions.\footnote{Sheaves already appear in Cartan's seminar of the year 1948. The theory was developed by Cartan and Oka, and in the year 1951, the famous  \emph{Cartan theorems A and B},\index{theorem!Cartan A and B} concerning a coherent sheaf on a Stein manifold, were  proved by Cartan, based on Oka's work. In 1954, Serre introduced sheaf theory in algebraic geometry, for varieties defined over algebraically closed fields.} Likewise, an affine variety is defined by its ring of ``regular" functions.  From the point of view of the theory of schemes\index{scheme}, one starts with an arbitrary commutative ring (with unit) $A$,\footnote{It may be useful to stress that it is after Grothendieck's ideas that commutative rings, instead of fields, played their central role in  algebraic geometry.} and one associates to it an affine scheme, namely, the set of prime ideals of $A$ equipped with a certain topology and with a sheaf structure whose fibers are the local rings defined by  these prime ideals. The main novelty with the Grothendieck setting is that classically the sheaf is a sheaf of functions. One may think of this in analogy with the manifolds defined by local charts, in the differential category.  A scheme\index{scheme} is a locally ringed space, modelled on an affine scheme. 
A locally ringed space\index{locally ringed space}\index{space!locally ringed} is a topological space equipped with a sheaf of commutative rings, called the structural sheaf of the space whose ring of germs at each point is a local ring. An affine scheme is a locally ringed space which is the spectrum of a commutative ring, that is, its set of prime ideals, equipped with the Zariski\index{Zariski topology}\index{topology!Zariski} topology.\footnote{The Zariski topology is a topology used in the study of algebraic varieties. The definitions in this topology are done in a purely algebraic way. We recall that Grothendieck introduced later on the \emph{\'etale topology},\index{etale topology}\index{topology!\'etale} a refinement of the Zariski topology.} It is equipped with a natural structural sheaf of rings which roughly correspond to germs of functions around elements. Schemes generalize algebraic varieties. Grothendieck reviews all these notions in Expos\'e II of his series of talks \cite{G2}; cf. \S \ref{sII} below. He introduced schemes  and locally ringed spaces by gluing together spectra of commutative rings using the Zariski topology. Grothendieck's development of these theories was partly motivated by his attempts to prove the Weil conjectures. We recall that Grothendieck's strategy was realized later on by Deligne, who also introduced new ideas.
 
Regarding also Item (\ref{44}), it is well known that the presence of nilpotent elements in structure rings of  schemes turned out to be an essential ingredient in algebraic geometry, and in particular in the theory of infinitesimal deformations of algebraic varieties. We recall that points in an analytic space are defined as ideals in rings. Grothendieck discovered that to  get the full strength of the theory, we must admit nilpotent elements in the structure rings.\index{nilpotent element} Thus, the equation $y^2=0$ defines a subscheme of the affine line with coordinate $y$ which is not the same as the subscheme defined by the equation $y=0$.
 To say it in simple words, a Riemann surface over the point $0$ is just a Riemann surface, but a Riemann surface defined over the equation $y^2=0$ is a Riemann surface together with a first-order deformation.\footnote{The fact of admitting nilpotent\index{nilpotent element} elements is stressed by Grothendieck in his talk at the Bourbaki seminar, at about the same period \cite{Gro-Bou182} p. 1: ``[...] In particular, whether or not  there is a given base field, there was no reason to exclude the case where these rings contain nilpotent elements. Until now, geometers refused to take into account these indications and persisted in restricting themselves to the consideration of affine algebras without nilpotent elements,\index{nilpotent element} i.e. of algebraic spaces whose structure sheaves do not contain any nilpotent element (and even, in most cases, to ``absolutely irreducible" algebraic spaces). The speaker considers that this state of mind was a serious obstacle to the development of truly natural methods in algebraic geometry." [[...] En particulier, qu'il y ait ou non un corps de base donn\'e, il n'y avait pas lieu d'exclure le cas o\`u ces anneaux contiennent des \'el\'ements nilpotents. Jusqu'\`a pr\'esent, les g\'eom\`etres s'\'etaient refus\'es \`a tenir compte de ces indications et se sont obstin\'es \`a se resterindre \`a la consid\'eration d'alg\`ebres affines sans \'el\'ements nilpotents, i.e. d'espaces alg\'ebriques dans les faisceaux structuraux desquels il n'y a pas d'\'el\'ements nilpotents (et m\^eme le plus souvent, des espaces alg\'ebriques ``absolument irr\'eductibles"). Le conf\'erencier pense que cet \'etat d'esprit a \'et\'e un obstacle s\'erieux au d\'eveloppement des m\'ethodes vraiment naturelles en g\'eom\'etrie alg\'ebrique].} In Grothendieck's theory of infinitesimal deformations,\index{infinitesimal deformation} nilpotent elements play a major role.\footnote{We note that in his infinitesimal calculus, Leibniz\index{Leibniz, Gottfried Wilhelm}  considered ``dual numbers," which contain the idea of nilpotent elements.} In general, the infinitesimal deformation theory of algebraic objects (algebraic curves, algebraic surfaces, etc.) defined over a field $k$ involves the study of objects defined over Artinian local rings with residue field $k$. Grothendieck's contribution is essential in this theory. 
 
 \medskip

 In his  construction of Teichm\"uller space, Grothendieck treats the case of closed surfaces of genus $\geq 2$ and he  leaves the construction of the Teichm\"uller space of genus $1$ as an exercise. Regarding the latter, he says that the reader can treat more generally the case of moduli spaces of complex tori or of complex abelian varieties, by a transcendental proof which is based on the fact that a complex torus can be regarded as a quotient of a vector space by a discrete subgroup of maximal rank, and a ``family of complex tori" above an analytic space $S$ as the quotient of a fiber bundle over $S$ associated to a locally free sheaf of rank $g$ by a discrete sub-sheaf of maximal rank. We shall say more about this in \S \ref{s:functor} below.

  In the rest of this survey, we shall review in detail some of the major ideas expressed in these lectures.

 \section{A review of Grothendieck's Expos\'e I: Teichm\"uller space representing a functor}\label{s:functor}
 
In this section, we review in some detail Grothendieck's first lecture \cite{G1}. In this lecture, he presents his goal, introduces the tools that he uses, and outlines the method of the proof of the existence and uniqueness of Teichm\"uller space representing a functor.

 It was usual for Grothendieck to give in advance the plan of something he was going to write several months (sometimes, years) later.  This was one of the characteristics of his broad and long-term vision and flux. This is the reason for which the outline of the proof that is given in this first lectures is not rigorously the one that he follows in the later lectures. In the meantime, he introduced new concepts and changed some steps in the general scheme of the proof. Therefore, the set of 10 written texts does not follow precisely the plan announced in the first lecture. In the last lecture \cite{G10}, Grothendieck gives the proof of the main results announced in the first lecture.  One must remember here that the Cartan seminar lecture notes were distributed gradually, week after week, and they were not altered a posteriori.
  
Let us make a few preliminary remarks on the proof;  some of them will be expanded in the rest of this text. 
 
One of the main tools is the introduction of categories whose objects are fibrations over complex spaces whose fibers are complex spaces. The base space may be fixed or variable, and depending on that, one gets different categories. Some of these categories are obtained as fiber products (or pull-backs) of other categories. The fibers may also be discrete (in this case the fibrations are covering spaces). There are several natural functors between these categories, with various properties. Some of the functors are \emph{rigidifying}\index{rigidifying functor}\index{functor!rigidifying} and some are \emph{representable}\index{representable functor}\index{functor!representable} and we shall explain these words. Rigidifying functors are obtained by equipping the fibers with some extra structure. For instance, when the fibers are Riemann surfaces, then equipping them with a marking is a rigidification of the structure. This extra structure is transported from fiber to fiber by analytic continuation, and in this process appears a monodromy, which is encoded (in the case in which we are interested) by the mapping class group action. The rigidified functor resulting from the rigidification of Riemann surfaces is the so-called \emph{Teichm\"uller functor}.\index{Teichm\"uller functor} This is an example of a representable functor. 
Grothendieck, in his series of lectures, proves the existence and uniqueness of the Teichm\"uller functor. The analytic space that represents it is a fiber bundle over Teichm\"uller space  which satisfies some universal property.

Let us now state in more precise terms Grothendieck's definitions and results. 

In what follows, all morphisms are complex analytic. 

\begin{definition} \label{def:over}
An \emph{analytic space $X$ over an analytic space $S$}\index{analytic space over an analytic space} is a pair $(X,f)$ where $X$ and $S$ are analytic spaces $X$ and $f:X\to S$ is an analytic map. The map $f$ is called the \emph{projection map}.
\end{definition}

We shall use (like Grothendieck) the notation $X/S$ for an analytic space over an analytic space $S$, highlighting the base $S$. 

Given an analytic space $X/S$, with projection map $f:X\to S$,  and a point $s$ in $S$, the \emph{fiber} over $s$ is the inverse image $f^{-1}(s)$. The fiber above any point is an analytic space. More precisely, it is the fiber product $X\times_S \{s\}$. (See below for the general definition.) The analytic space $X/S$ over $S$ is sometimes denoted by $(X,f)$. The pair $(X,f)$ is regarded as a family of analytic spaces parametrized by $S$, and sometimes, a \emph{deformation} of one of the fibers of the map $f$. 

In the case at hand, the fibers are Riemann surfaces of the same genus. Under the technical condition that the map $f$ is flat, we may think of this map as a deformation of the complex structure of some fiber. 

Families of analytic spaces of dimension $2$ or more are usually studied within the theory of Kodaira-Spencer. Grothendieck also developed a theory that included this deformation theory into the abstract algebraic geometry setting. Fiber products and the related theory of ``base change," in algebraic topology, are central objects in Grothendieck's theory.

Given an analytic space $X/S$ and an open set $U$ in $S$, the \emph{restriction}\index{analytic space over an analytic space!restriction} $(X\vert_U)/U$ of $X/S$ to $U$ is the analytic space $f^{-1}(U)$ equipped with its natural fibration over $U$.

For a fixed analytic space  $S$, the set of analytic spaces over $S$ form a category $\mathcal{A}_S$. A \emph{morphism} between the two objects $(X,f)$ and $(X',f')$ is an analytic map $g:X\to X'$ commuting with the projections and the identity map on $S$. In other words, we have $f'\circ g=f$, or the following diagram commutes:

 \[\begin{array}[c]{ccc}
 X & \stackrel{g}{\longrightarrow}& X'\\
 \scriptstyle{f}  \downarrow && \downarrow \scriptstyle{f'}\\
  S & \stackrel{\mathrm{Id}}{\longrightarrow}& S'\\
\end{array}\]

\medskip

For any morphism between analytic spaces \[h:S'\to S,\]
we have a functor 
\[\mathcal{A}_S\to \mathcal{A}_{S'}.\] 

Grothendieck considers this morphism as a \emph{base-change functor}. It sends every analytic space
 $f:X\to S$ over $S$ to the analytic space $X\times_S S'\to S'$ over $S'$ where $X\times_S S'$ is the \emph{fiber product}\index{fiber product} of $X$ and $S'$ over $S$. We recall that as an underlying topological space, this is the set of pairs $(x,y)\in X\times S'$ satisfying $f(x)=h(y)$,  where the map $X\times_S S'\to S'$ is induced from the projection on the second factor. The introduction of the sheaf defining the analytic structure of this space needs some care. The fiber product $X\times_S S'\to S'$ is also called the \emph{pull-back} \index{analytic space over an analytic space!pull-back}\index{pull-back} of $X$ by $h$ and denoted by $h^* X$. It is an analytic space over $S'$. 

The notion of fiber product is a generalization of the notion of extension for scalars.

An analytic space $X$ over $S$ is said to be \emph{proper}\index{analytic space over an analytic space!proper} if the associated map $f:X\to S$ is proper, that is, if the inverse image by $f$ of any compact subset of $S$ is a compact subset of $X$. Note that this is a topological condition which does not use the analytic structures of the spaces $X$ and $S$. The space $X$  is said to be \emph{simple} if for every $x$ in $X$ there exists a neighborhood $U$ of $s=f(x)$ and a neighborhood $V$ of $x$ over $U$ such that $V$ is $U$-isomorphic to a product $U\times C$ where $C$ is a complex variety. This condition easily implies that there are no other singularities on the fibers than those that arise from $U$. The last statement is a kind of an implicit function theorem \`a la Grothendieck.

\begin{definition} An \emph{analytic curve of genus $g$} over $S$ is an analytic space $X$ over $S$ which is proper and simple and such that the fibers (which are compact complex varieties) are connected Riemann surfaces of genus $g$.
\end{definition}

\begin{lemma} \label{lem:G} An analytic curve of genus $g$ over an analytic space $S$ is a locally trivial topological fiber space.
\end{lemma}

This lemma is easy to prove if $S$ is a manifold,\footnote{A theorem attributed to Ehresmann says that if $\pi:X\to S$ is a proper submersion between differentiable manifolds, then $\pi$ is locally trivial in the $C^{\infty}$ sense. For our purposes, it is not sufficient to consider the case where $S$ is a manifold, and one must also deal with singular complex spaces.} but in the present setting, $S$ might not be a manifold.  The lemma is proved by Grothendieck in a later lecture (\cite{G10}, Proposition I.8).

With the above definitions, the main question becomes that of understanding the space of isomorphism classes of analytic curves over a given  space $S$.

To do this, Grothendieck introduces several categories. The first of them is the category  $\mathcal{F}_S$ of algebraic curves of genus $g$ over a fixed space $S$.
The second category, denoted by $\mathcal{F}$, is the category of algebraic curves of genus $g$ over variable spaces. This is a \emph{fibered category},\index{fibered category}\index{category!fibered}  that is, a category in which the operation of fiber product is defined.\footnote{The reader may remember that the set-theoretical fiber product is always defined, but that we ask here that the resulting objects are elements of the category, that is, they are analytic, etc. With this requirement, such an operation does not exist in an arbitray category.} In Grothendieck's words, the aim of his seminar talks is ``to determine, in a simple way, the structure of the fibered category $\mathcal{F}$, up to fibered category equivalence."  In fact, the problem considered by Grothendieck amounts to that of finding a \emph{universal} object, also called a \emph{final} element in this second category. Roughly speaking, this is an object  from which every other object is obtained by pull-back.\footnote{For this to hold, the morphisms of $\mathcal{F}$ should be defined with care; a morphism from $X\to S$ to $X'\to S'$ is given by a morphism from $S$ to $S'$ and an \emph{isomorphism} between $X$ and $X' \times _{S'} S$ (that is, $\mathcal{F}$ is a category which is fibered in groupoids).} We shall be more precise below.

Let $\mathcal{A}$ be the set of objects of $\mathcal{F}$ up to isomorphism and $\mathcal{A}(S)$ the set of objects of $\mathcal{F}_S$ up to isomorphism. In other words, $\mathcal{A}(S)$ is the set of curves of genus $g$ up to isomorphism over a fixed space $S$.  

The dependence on $S$ of analytic curves over $S$ is a contravariant functor 
\[\mathcal{A}: \mathbf{(An)}\to  \mathbf{(Ens)}\]
from the category $\mathbf{(An)}$ into the category $\mathbf{(Ens)}$: to the analytic space $S$ we associate the set of analytic curves above $S$, and to a morphism $h:f: S'\to S$ between analytic spaces we associate the set-theoretic map which assigns to an analytic curve over $S$ its pull-back by $f$ on $S'$. This pull-back is an analytic curve over $S'$.

The first approach to the solution of the moduli problem is to obtain a \emph{representation} of the functor $\mathcal{A}$. Grothendieck notes that if no object of $\mathcal{F}_S$ had a nontrivial automorphism, then the category $\mathcal{F}_S$  would be known, up to equivalence, by the set $\mathcal{A}(S)$. In other words, if no object of $\mathcal{F}_S$ has non-trivial automorphisms, then the map $\mathcal{F}_S\to \mathcal{A}(S)$ is injective, and the morphisms in $\mathcal{F}_S$ are also determined by this map.  Thus, the category $\mathcal{F}_S$  would be known by the set represented. 
 And if the same  condition were satisfied for every $S$, then the fibered category $\mathcal{F}$ would be known up to equivalence from the functor $\mathcal{A}$. If furthermore the functor $\mathcal{A}$ were \emph{representable}, that is, if it were isomorphic to a functor of the form 
\[S\to  \mathrm{Hom}(S,M)\] from the category $\mathbf{(An)}$ of analytic spaces into the category $\mathbf{(Ens)}$ of sets,
where $M$ is some appropriate object of $\mathbf{(An)}$ ($M$ will be in fact Teichm\"uller space), then the fibered category $\mathcal{F}$ would be known up to equivalence from that object $M$, because everything would be obtained by a pull-back, using a map onto $M$.  

It turns out, as Grothendieck notes right at the beginning (p. 7-03) that the functor $\mathcal{A}$ is not representable, and the reason is that some curves of genus $g$ may have non-trivial automorphisms. The  goal of the theory is then to eliminate these automorphisms. This leads to the introduction of another functor, denoted by $\mathcal{P}$. This is a representable functor from the category of marked curves of genus $g$ over analytic spaces to the category of sets, and the object that represents this category is the Teichm\"uller curve. This is the content of Grothendieck's theorem. The result is equivalent to Teichm\"uller's theorem (Theorem \ref{th:T}) which says that one can reconstruct any globally analytic family over an analytic space $\mathfrak{B}$ by taking the pull-back by a map $f:\mathfrak{B}\to\mathfrak{C}$ of the Teichm\"uller curve $\underline{\mathfrak{H}}[\frak{c}]$ over $\mathfrak{C}$.

With the existence of the functor $\mathcal{P}$, the categories $\mathcal{F}$ and $\mathcal{F}_S$ become \emph{representable fibered categories}. 

The mapping class group appears now in the construction of principal coverings of $S$.
We start by describing a general topological construction which Grothendieck presents;  it uses Lemma \ref{lem:G}.

Let $X$ be a topological space equipped with a locally trivial fibration over a topological space $S$, a fiber being a finite simplicial complex $C_0$, and let $G$ be the extended mapping class group of $C_0$, that is, the group of homotopy classes of homeomorphisms of this space. 

For a given point $s$ in $S$, let $X_s$ be the fiber above it. The space $I(C_0,X_s)$ is the set of homotopy classes of homeomorphisms between the two spaces. It is a principal homogeneous space for the group $G$. The union of the sets $I(C_0,X_s)$, for $s$ in $S$, has a natural structure of a principal covering above $S$, with group $G$. It is denoted by $\mathcal{R}(X/S)$. For a fixed $C_0$, Grothendieck notes the following properties of the operation $\mathcal{R}$:
\begin{enumerate}
\item For $S$ fixed, this is a functor in $X$.
\item  It is compatible with the operation of taking pull-backs of fiber spaces.
\end{enumerate}

The group $G$ also acts on the homology and the cohomology groups of $C_0$, and one obtains the local systems of the homology and the cohomology of the fibers of $X/S$ by taking the fiber bundles associated to the principal covering $\mathcal{R}(X/S)$, with the various operations described.

Now we specialize to the case where $C_0$ is a topological surface which underlies the algebraic curves of genus $g$ that we are considering.  The mapping class group of $C_0$ is the group of isotopy classes of orientation-preserving diffeomorphisms of this curve, and it is denoted by $\gamma$. It is an index-2 subgroup of the group $G$ associated to $C_0$ that were defined above. Grothendieck calls $\gamma$ the  \emph{Teichm\"uller group}.\index{Teichm\"uller group}\index{mapping class group}

Choosing continuously a system of orientations on the fibers $X_s$ is equivalent to choosing a section of the principal $\mathbb{Z}/2\mathbb{Z}$-covering of $S$ associated to $\mathcal{R}(X/S)$, that is, restricting the structure group $G$ to $\gamma$. The result is a principal covering of $S$, whose group is $\gamma$, and which is denoted by $\mathcal{P}(X/S)$.

As already noted, because curves of genus $g\geq 2$ may have non-trivial automorphisms, the functor $\mathcal{F}$ is not representable.  Furthermore the theory developed later shows that the existence of such automorphisms also prevents the functor $\mathcal{A}$ of being representable.

We now recall that in the particular case where the base $S$ is a point (that is, when the fibered spaces are ordinary curves of genus $g$), Teichm\"uller solved the problem of the elimination of nontrivial automorphisms of Riemann surfaces by introducing a marking,\index{marking!by collections of curves} and he also introduced the mapping class group action on the space of marked surfaces. The markings that he used for Riemann surfaces are not always mapping class group elements; he also used be decompositions of the surface by collections of simple closed curves such that an element of the mapping class group which sends each curve to itself up to homotopy is necessarily the identity mapping class. The non-trivial automorphisms of the structures considered are the cause of singularities. The elimination of these singularities is also done through the  marking;\index{marking!rigidification} in this sense the marking is a rigidification\index{rigidification} of the objects by adding an extra structure.  The mapping class group action becomes part of the structure. The introduction of this group then leads to the construction of an infinite covering space of moduli space, which is Teichm\"uller space. This is an example of a general idea of a ``rigidification by taking a covering space." It leads to the following definition due to Grothendieck:

The Teichm\"uller covering $\mathcal{P}(X/S)$ of $X/S$ is a  topological covering of $S$ and the mapping class group $\gamma$ acts on the fibers by interchanging sheets over points. This action is discrete and $\mathcal{P}(X/S)$ is a principal bundle over $S$. The covering $\mathcal{P}(X/S)$ is also equipped with an analytic structure, defined by the property that the local projections are analytic isomorphisms. The group $\gamma$ operates analytically on $\mathcal{P}(X/S)$ by $S$-automorphisms (that is, automorphisms that induce the identity on $S$).

The map 
\[X/S\mapsto \mathcal{P}(X/S)\]
is a covariant functor between the categories $\mathcal{F}_S$ and $\mathbf{(Ens)}$.  It associates to a morphism $h:X/S\to X'/S$ a map between coverings $\mathcal{P}(X/S)\to \mathcal{P}(X'/S)$ obtained by transporting by the homotopy class of $h$ a mapping class group action on the fibers of $X/S$ to a mapping class group action on the fibers of $X'/S$.

The passage from $X/S$ to $\mathcal{P}(X/S)$ has the following rigidity property which Grothendieck proves in his later talks:

\begin{theorem} \label{th:rigid}
Any automorphism of an $X/S\in \mathcal{F}_S$ which induces the identity on $\mathcal{P}(X/S)$ is the identity.
\end{theorem}

An equivalent form of Theorem \ref{th:rigid} is the following which we shall see in more detail below:

\begin{theorem} \label{th:rigidifying}
The functor $\mathcal{P}$ is a \emph{rigidifying functor}. 
\end{theorem}

We recall that this means that\index{rigidifying functor}\index{functor!rigidifying} every automorphism of an object of $\mathcal{F}_S$ which induces the identity on $\mathcal{P}(X/S)$ is the identity.

\begin{definition} A \emph{Teichm\"uller structure} (or a \emph{$\mathcal{P}$-structure}) on a curve $X$ of genus $g$ above an analytic space $S$ is a section of the Teichm\"uller covering $\mathcal{P}(X/S)$ over $S$.
 
\end{definition}
 Since $\gamma$ operates on the right on $\mathcal{P}(X/S)$, it also operates on the right on the set of Teichm\"uller structures.

\begin{definition}[$\mathcal{P}$-curve]
A \emph{$\mathcal{P}$-curve} is a curve of genus $g$ equipped with a $\mathcal{P}$-structure. 
\end{definition}

The $\mathcal{P}$-curves of genus $g$ form a fibered category, denoted by $\mathcal{F}^{\mathcal{P}}$, which lies above the category of analytic spaces. A \emph{$\mathcal{P}$-curve} above $S$ (that is, an element of $\mathcal{F}^{\mathcal{P}}$) does not have any nontrivial automorphism.

If $S$ is connected and nonempty and if the set of Teichm\"uller structures on $X/S$ is nonempty, then this set is a principal homogeneous set with structural group $\gamma$. 

\begin{definition} 
A curve of genus $g$ over $S$, equipped with a Teichm\"uller structure, is called a \emph{Teichm\"uller curve of genus $g$ over $S$.}
\end{definition}

From Theorem \ref{th:rigid}, we obtain:
\begin{proposition} An automorphism of a Teichm\"uller curve is the identity.
\end{proposition}

Using the notion of pull-back of a Teichm\"uller curve $X/S$ by a morphism $S'\to S$, we have:
\begin{proposition}
Teichm\"uller curves of genus $g$ form a fibered category over the category $\mathbf{(An)}$ of analytic spaces. 
\end{proposition}

Grothendieck states several properties valid in a much more  general setting. 

Let $\gamma$ be now a discrete group. We assume that for every analytic space $S$ we have a covariant functor 
\[\mathcal{P}: X/S\to \mathcal{P}(X/S)\]
from the category $\mathcal{F}_S$ of curves of genus $g$ above $S$ to the category of principal coverings of $S$ of group $\gamma$.
 
 The curves are now equipped with a marking such that every automorphism of the curve over $S$ which preserves the marking is the identity.

An automorphism of an object in $\mathcal{F}_S$ induces an automorphism on its image $\mathcal{P}(X/S)$.

In the case where $S$ is a point, the functor associates to a curve  a discrete space on which the group $\gamma$ acts transitively.  

Grothendieck gives the following examples of rigidifying functors.

Suppose we have two groups $\gamma$,  $\gamma'$ and a homomorphism $\gamma \to \gamma'$, and suppose we have a rigidifying functor  that uses the group $\gamma$, in the above sense.
If the kernel $\gamma\to \gamma'$ is small enough, then $\mathcal{P}'$ is also a rigidifying functor.

We get a functor called  \emph{fibering functor}\index{functor!fibering}\index{fibering functor} (``foncteur fibrant") from the fibered category of curves of genus $g$ into the fibered category of principal coverings of group $\gamma'$,  by defining:

\[\mathcal{P}'(X/S) = \mathcal{P}(X/S)\times_{\gamma}\gamma'.\]
This functor amounts to adding more marking.

Grothendieck shows that any rigidifying functor can be obtained in the above way from the Teichm\"uller functor.

In the next example, $\mathcal{P}$ is Teichm\"uller's rigidifying functor, and $\gamma$ is again the mapping class group.  

The group $\gamma$ acts on the homology group $H_1(C_0,\mathbb{Z})$ where $C_0$ is a fiber, leaving invariant the symplectic intersection form. Choosing a basis, we get a representation
\[\gamma \to \mathrm{Sp}(2g,\mathbb{Z}).\]

We get a fibering functor :
\[X/S\mapsto \mathcal{P}(X/S)\times_{\gamma} \mathrm{Sp}(2g,\mathbb{Z})= \mathcal{P}'(X/S).\]

If we replace the coefficient group $\mathbb{Z}$ by $\mathbb{Z}/n\mathbb{Z}$, we get a fibered functor with respect to the finite group $\mathrm{Sp}(2g, \mathbb{Z}/n\mathbb{Z})$:

\[X/S\mapsto \mathcal{P}(X/S)\times \mathrm{Sp}(2g,\mathbb{Z}/n\mathbb{Z})=\mathcal{P}_n(X/S).\]

Sections of this fiber space over $S$ are Riemann surfaces together with the choice of a symplectic basis of the homology of the surface with coefficients in $\mathbb{Z}/n\mathbb{Z}$. 
\begin{lemma}
 There exists an integer $n>0$ such that the fibered functor $\mathcal{P}_n$ is rigidifying.
\end{lemma}

This implies that $\mathcal{P}'$ is rigidifying, and it then implies that  $\mathcal{P}$ is rigidifying. This is used in the proof of Theorem \ref{the:GG}.

 $\mathcal{F}_S^{\mathcal{P}}$ is \emph{the functor $\mathcal{F}_S$ rigidified by $\mathcal{P}$}. 
For any analytic space $S$, let $\mathcal{A}_{\mathcal{P}}(S)$ be the set of objects of $\mathcal{F}^{\mathcal{P}}_S$ up to isomorphism.
Then,
\begin{theorem}The contravariant functor
\[S\mapsto \mathcal{A}_{\mathcal{P}}(S)\]
is representable.
\end{theorem}

  \section{A quick survey of Grothendieck's Expos\'es II to IX}\label{sII}
 
In this section, we briefly describe the content of Expos\'es II to IX of Grothendieck's lectures at Cartan's seminar. 

\bigskip

 \noindent{\bf Expos\'e II. G\'en\'eralit\'es sur les espaces annel\'es et les espaces analytiques.} (General facts on ringed spaces and analytic spaces.) \cite{G2}

The goal in this talk is to introduce in a natural way the notion of analytic space. The notion of analytic manifold was already introduced by Cartan and others, and Grothendieck needs to introduce the more general notion of analytic variety, which may be singular.  This is needed because the moduli space $M_g$  of curves of genus $g$ is not a smooth manifold. (In fact, it is not even a manifold.)
In this expos\'e, Grothendieck explains how to allow nilpotent elements in the structural sheaf, on the model of what he did for schemes.

 Before considering analytic spaces, Grothendieck introduces the more general notion of ringed space\index{ringed space} (``espace annel\'e").  The main idea is that points are defined as local functions.

   A \emph{ringed space}\index{ringed space}\index{space!ringed} is a topological space\footnote{The topology on $X$ may be an unusual topology, for instance the Zariski topology of an algebraic variety, which is generally not Hausdorff, or the etale topology (introduced by Grothendieck).} $X$ equipped with a collection of commutative rings $\mathcal{O}_X$ that forms the \emph{structure sheaf} of $X$. The elements of these commutative rings can be thought of as functions defined on open sets of $X$. In a ringed space, a point is described as a maximal ideal. Thus, a ringed space has an underlying topological space but it carries much more information than the topology of this underlying space.
   
   A local ring is a commutative ring which has a unique maximal ideal. A typical example is the ring of rational functions on $\mathbb{C}$ for which a given point is not a pole. This is the ``local ring of rational functions" at the given point. Grothendieck then introduces  the notion of a \emph{locally ringed space}\index{locally ringed space}\index{ringed space!locally} (``espace annul\'e par anneaux locaux") and of (local) morphisms between such spaces. A ringed space $(X,\mathcal{O}_X)$ is locally ringed when the fibers of the structure sheaf $\mathcal{O}_X$ are local rings.   
   
   An \emph{open subset} of a locally ringed space $X$ is an open subset $U$ of the underlying topological space equipped with the restricted sheaf of rings $\mathcal{O}_X\vert_U$. The pair $(U,\mathcal{O}_X\vert_U)$ is then itself a locally ringed space. 
   
An example of a locally ringed space is a topological space $X$ with structure sheaf the sheaf $\mathcal{O}_X$ of real-valued or complex-valued continuous functions on open subsets of $X$. Other examples include differentiable manifolds, complex manifolds and analytic varieties.

 The spectrum\index{spectrum!commutative ring} of a commutative ring $A$ (that is, the set of its prime ideals) is a locally ringed space. This space is equipped with a topology obtained by declaring the closed sets to be the collections of prime ideals of $A$ that contain a given ideal. We recall that a scheme\index{scheme} is a locally ringed space obtained by gluing together spectra of commutative rings.  The topological space $X$ defined on the spectrum of a commutative ring has a naturally defined sheaf of rings, its structure sheaf $\mathcal{O}_X$, and the ringed space   $(X,\mathcal{O}_X)$ is called an \emph{affine scheme}. The category of affine schemes is equivalent to the (dual) category of commutative rings. One of the goals of Grothendieck was to build a geometric category which is equivalent to  the (dual) category of commutative rings.
      
         There is a natural notion of morphism between ringed spaces or locally ringed spaces:   If $(X,\mathcal{O}_X)$ and $(Y,\mathcal{O}_Y)$ are two locally ringed spaces, a morphism between them is given by a continuous map $f:X\to Y$ together with a morphism $f^\sharp: \mathcal{O}_Y\to f^* \mathcal{O}_X$ from the structure sheaf of $Y$ to the direct image of the structure sheaf of $X$ such that for every $x$ in $X$, the ring morphism $\mathcal{O}_{Y, f(x)}\to \mathcal{O}_{X,x}$ induced by $f^\sharp$ is a morphism of local rings, that is, it sends the maximal ideal of $\mathcal{O}_{Y, f(x)}$ to the maximal ideal of $\mathcal{O}_{X,x}$. Thus, a morphism between (locally) ringed spaces is specified by a pair of maps $(f,f^{\sharp})$. Equipped with these morphisms, ringed spaces and locally ringed spaces form two categories, the second one being a subcategory of the first.

 In the analytic case,\footnote{In the algebraic case, the definiton has to be modified.} a morphism  $f: (X,\mathcal{O}_X)\to (Y,\mathcal{O}_Y)$ between locally ringed spaces is said to be unramified\footnote{Algebraic geometers say ``unramified" (and they also talk of a ``net morphism") for something differential topologists would call an immersion.} if $f$ is a topological immersion and if for all $x$ in $X$ the ring morphism $\mathcal{O}_{Y, f(x)}\to \mathcal{O}_{X,x}$ is surjective.
   
A morphism between ringed spaces is not determined by the induced map on the underlying topological spaces, and as always in this setting, a morphism in which the underlying topological spaces are points is already an interesting object. For instance, for every point $x$ in a locally ringed space $(X,\mathcal{O}_X)$, the topological space $\{x\}$ equipped with the constant sheaf $k(x)$ is a  locally ringed space. We have an immersion between locally ringed spaces $(x,k(x))\to (X,\mathcal{O}_X)$.

  Examples of locally ringed spaces include:
  
\begin{itemize}

\item Differentiable manifolds, equipped with their $\mathbb{C}$-valued differentiable functions;

\item Complex varieties, equipped with their holomorphic functions;
    
\item Algebraic varieties;
\item Schemes, which we already encountered. 
\end{itemize}

More precisely, a scheme\index{scheme} is a locally ringed space $(X,\mathcal{O}_X)$ such that any point $P$ of $X$ has a neighborhood $U$ such that the restricted locally ringed space $(U,(\mathcal{O})_{\vert U})$ is an affine scheme.

An important idea in Grothendieck's theory of schemes is that a scheme is characterized by the morphisms from other schemes into it.

   There is a tangent space associated to a locally ringed space which generalizes the standard construction for differentiable manifolds. Infinitesimal deformations correspond to enlarging the structural ring.

Let us now review the content of Expos\'e II.
 
 Grothendieck declares in the introduction to Expos\'e II  that his aim in the series of talks is to give a sketch, as brief as possible, of the elements of the foundations of analytic geometry that are needed in the formulation of the most important existence problems, with the aim of proving the theorems stated in Expos\'e I. The theorems generalize the existence and uniqueness theorem of Teichm\"uller (Theorem \ref{th:T} above). Except for some special cases, the results are developed in a setting which is more general than that of complex algebraic geometry, namely, they are valid for algebraic geometry over an arbitrary complete field equipped with a valuation. This will play an important role in Grothendieck's general program
 on moduli spaces of varieties over general fields. The special cases are those of Expos\'e VII \cite{G7}, in which Grothendieck presents what he considers as ``more profound results," namely, on the theory (of the type GAGA)  of projective morphisms, due to Grauert-Remmert, and the finiteness theorem of Grauert, which are special to the case where the base field is the field $\mathbb C$ of complex numbers.

 Some pre-requisites are needed in the definition of an analytic space. Quoting Grothendieck: 
\begin{quote}\small
Our immediate goal is to develop certain global existence theorems and at the same time the language in which they can be formulated. It appears that the list of result which we shall need amounts to the following (\cite{Cartan1} \cite{Cartan2}):

\begin{enumerate}
\item The ring $A_n=k(t_1,\ldots,t_n)$ of convergent power series in $n$ variables is a Noetherian local ring whose completion is isomorphic to the ring of formal power series in the $t_i$.

\item Let $A\to B$ be a homomorphism of $k$-analytic local algebras over $k$ (that is, isomorphic to non-zero quotients of the rings $A_i$), such that the competed ring $\hat{B}$ is a finite-type module over the completed ring $\hat{A}$ (we say that $B$ is quasi-finite over $A$). Then $B$ is a finite-type module over $A$.

\item The sheaf of analytic functions of the space $k^n$ is a coherent ring sheaf.
\end{enumerate}
\end{quote}

The result in (1) replaces the implicit function theorem in this non-smooth setting. The inverse may not exist in the setting of polynomials, and one uses instead this notion of completion of the rings.

Property (2) is a lifting property.

property (3) is a coherence property.  The notion of coherence is close to the notion of being of finite type.

These results are used to introduce the following notion of analytic space. 

Given a field $k$ equipped with a complete valuation, a positive integer $n$ and given an open set of $k^n$, an \emph{analytic function} on $U$ is a map $U\to k$ defined by a convergent power series. When $U$ varies we get a sub-sheaf from the sheaf of maps from $k^n$ to $k$, called the \emph{sheaf of holomorphic functions} on $k^n$. Equipped with this sheaf, $k^n$ becomes a ringed space denoted by $\mathbb{E}^n$. 

Then Grothendieck makes the following definition:

\begin{definition}[Definition 2.1 of \cite{G2}] An analytic space over a complete field $k$ equipped with a valuation is a topological space  ringed  by $k$-algebras, such that every point has a neighborhood which is isomorphic, with respect to the induced structure, to a finitely presented ringed sub-space of a space $\mathbb{E}^n$. 
\end{definition}
 
 There are natural definitions of analytic subspaces and of morphisms between them. Analytic spaces form a category. The first examples of analytic spaces are the space $\mathbb{E}^n$ defined above and analytic manifolds. Grothendieck gives several other examples from algebraic geometry.

 \medskip
 \noindent{\bf Expos\'e III. Produits fibr\'es d'espaces analytiques.} (Fibered products of analytic spaces.) \cite{G3}

In this Expos\'e  Grothendieck introduces a functor $\phi$ 
\[\phi(X) : \mathrm{Hom}(X,\mathbb{E}^n)\to \Gamma(X,\mathcal{O}_X)^n, \]
defined by composing with the coordinate functions. Here, $\mathbb{E}^n$ is the sheaf of holomorphic functions; this is  the ``model space" defined in Expos\'e II.  

Grothendieck calls this functor the \emph{sheaf homomorphism}, and he proves that it is an isomorphism (Theorem \ref{th:sheaf}). He then determines  the morphisms from an arbitrary analytic space $X$ to $\mathbb{E}^n$.
The functor $\phi$ is needed to map curves in projective space. This is related to Hilbert schemes.

 For $n\geq 1$, the set $\Gamma(X,\mathcal{O}_X)^n$ is the set of $n$ sections of the structure sheaf $\mathcal{O}_X$ of $X$.  For $n=1$, this set is simply the set of functions on $X$.
Grothendieck proves the following.
\begin{theorem}[Theorem 1.1 of \cite{G3}] \label{th:sheaf}
For any analytic space $X$, the sheaf homomorphism $\phi$
is an isomorphism.
\end{theorem}

\begin{corollary}[Corollary 1.2 of \cite{G3}] The contravariant functor
\[
X \to \Gamma(X,\mathcal{O}_X)^n
\]
defined on the category $\mathbf{(An)}$ of analytic spaces with values in the category $\mathbf{(Ens)}$ of sets is representable by the object $\mathbb{E}^n$. 
\end{corollary}

\begin{theorem}[Theorem 2.1 of \cite{G3}] 
In the category $\mathbf{(An)}$ of analytic spaces, finite projective limits exist. The canonical functor 
\[T: \mathbf{(An)}\to \mathrm{(Top)},\] which associates to every analytic space the underlying topological space, commutes to projective limits.
 \end{theorem}

Grothendieck then proves existence theorems for fiber products in the analytic category,  that is, the product of  two spaces equipped with morphisms hitting the same target. These results are needed for base change, and this is important for comparing families, and also for passing to finite covers, and undoing monodromies. 

The infinitesimal structure (that is, the local ring structure) on the fiber product is obtained from the infinitesimal structure on each of the two factors. 
This is also useful in gluing local deformation spaces.

 \medskip
 
 \noindent{\bf Expos\'e IV. Formalisme g\'en\'eral des foncteurs repr\'esentables} (General formalism of representable functors.) \cite{G4}
  We already outlined the ideas in this section in \S \ref{s:ideas} of the present chapter. 
 
  \medskip
 \noindent{{\bf Expos\'e  V. Fibr\'es vectoriels, fibr\'es projectifs, fibr\'es en drapeaux.} (Vector bundles, projective bundles, flag bundles.) \cite{G5}

A flag bundle\index{flag bundle} is defined as an increasing sequence of sub-bundles, each one being locally a direct factor of the next one, generalizing the usual notion of flag in linear
algebra to the setting of ringed spaces. This leads to the theory of a ``relative scheme"\index{relative scheme} over a ringed space.  The algebraic and topological operations and the algebraic and topological languages enter into the realm of analytic geometry. Grothendieck  explains how to associate to vector bundles, projective bundles and flag bundle functors which are representable. 

The paper starts with some preliminary remarks saying that it is unfortunate that the language of  schemes relative to a ringed space was not developed in the foundational essay \cite{Dieu-GroI} by Dieudonn\'e and Grothendieck, that one has to work in the category of analytic spaces but that there is a clear necessity to develop a relative\footnote{A relative theory is another name for a theory of spaces over spaces. We already mentioned examples.} theory of analytic spaces over general spaces ringed by topological rings. This is an  instance of Grothendieck's idea of establishing functors between theories so that constructions and results in one theory follow from constructions and results in the other one. In the present setting, the ``relative theory" would be used in the following way:

To each algebraic space $X$ and to each relative scheme of finite type over $S$, one can associate an analytic space $X_{\mathrm{an}}$ such that for any analytic space over $S$ (that is, an analytic space equipped with an analytic morphism $X\to S$; see definition \ref{def:over} below) one has a bijection $\mathrm{Hom}_S(T, X_{\mathrm{an}}) \rightarrow \mathrm{Hom}_S(T,X)$ which is functorial in $T$. Here, the first Hom is understood in the category $\mathbf{(An)}$ of analytic spaces and the second one in the category of ringed spaces  with local rings. Grothendieck declares that ``all the operations of vector bundles, flag bundles, etc. in analytic geometry can be deduced from analogous constructions in algebraic geometry by applying the functor $X\to X_{\mathrm{an}}$, and their main formal properties trivially follow from results in algebraic geometry and from the elementary properties of the functor $X\to X_{\mathrm{an}}$."

 \medskip
 \noindent{{\bf Expos\'e VI. \'Etude locale des morphismes : germes d'espaces analytiques, platitude, morphismes simples.} (Local study of morphisms: germs of analytic spaces, flatness and simple morphisms.) \cite{G6}

Grothendieck gives the definition of a germ of an  analytic space and of morphisms between two such spaces. He then defines a functor 
\[\mathrm{(germes \ An})^{\mathrm{o}} \to \mathrm{(Alg \ An)}
\]
from the category of germs of analytic spaces into the category of analytic algebras
and he proves that this gives an equivalence between the inverse category of the category of germs of analytic spaces and the category of analytic algebras.

 \medskip
 \noindent{{\bf Expos\'e VII. \'Etude locale des morphismes: \'el\'ements de calcul infinit\'esimal.} (Local study of morphisms: elements of infinitesimal calculus.) \cite{G7}
 
 This expos\'e concerns local deformations.

Grothendieck develops the use of ringed spaces with nilpotent elements. He formulates fundamental geometric notions and the elements of the infinitesimal geometry of analytic spaces over analytic spaces with nilpotent elements that appear in the geometric constructions that will follow.
 In particular, he defines the notion of submersion in this non-smooth case, in which infinitesimal neighborhoods are described by powers of the defining ideals. He also gives a differential characterization of local immersions and he provides several examples.

 \medskip
 \noindent{{\bf Expos\'e  VIII. Rapport sur les th\'eor\`emes de finitude de Grauert et Remmert.} (Report on the finiteness theorems of Grauert and Remmert.) \cite{G7}

The results in this talk, unlike those of the rest of the series, are restricted to the case where the ground field is the field $\mathbb{C}$ of complex numbers.
Here, Grothendieck presents results of Grauert \cite{Grauert} and Grauert-Remmert \cite{Grauert-Remmert} on morphisms between analytic spaces. These results belong to a series of ``finiteness theorems" in analytic geometry. Grauert's finiteness result \cite{Grauert} says that a certain module associated to proper morphisms between analytic spaces is coherent. The coherence property guarantees that certain cohomology groups are finite-dimensional. One then deduces the desired finiteness property. Grauert's finiteness result is a generalization of the so-called Cartan-Serre finiteness theorem \cite{Cartan-Serre} and \cite{Serre} of the type ``GAGA," stating that the cohomology vector spaces of a compact complex space with values in a coherent sheaf are finite-dimensional.

  Grothendieck declares that this theorem should be interpreted as an expression of the equivalence between the ``algebraic geometry of relative projective schemes over an analytic space $Y$" and the theory of ``projective analytic spaces over $Y$." He also says that the analogues of these two theorems in the setting of algebraic geometry were already known and much more easier to prove. Grothendieck refers to \cite{Dieu-GroIII}.

We note that there is a Bourbaki seminar by Cartan on Grauert's results \cite{Cartan-Grauert}.

 \medskip
 \noindent{{\bf Expos\'e  IX. Quelques probl\`emes de modules.} (Some moduli problems.) \cite{G9}

In this penultimate talk, Grothendieck reviews some typical ``moduli" problems. These include Hilbert moduli and Picard moduli.

For the existence theorem of Hilbert moduli spaces in the setting of schemes, he refers to \cite{Gro-BouIV}. This theorem is the essential tool that is needed in Expos\'e X for the construction of Teichm\"uller space.

\section{A review of Grothendieck's Expos\'e X: The construction of moduli space and of Teichm\"uller space}

  We now review Grothendieck's last lecture in the series, \cite{G10}, in which he describes the construction of Teichm\"uller space and gives his own version of Teichm\"uller's theorem (Theorem \ref{th:T}). The proof uses what Grothendieck calls the ``rigidity of the Jacobi functor of level $n\geq 3$ ." The statement of this last result is the following:
\begin{theorem}[Theorem 3.1 of \cite{G10}] For any integer $n\geq 3$, the Cartesian Jacobi functor of level $n$ defined on the fibered category of curves of genus $g$ on analytic spaces $S$ is a rigidifying functor.
\end{theorem}

The Jacobian\index{Jacobian} $J(S)$ of the closed surface $S$ is the torus $ H_1(S,\mathbb{R})/H_1(S,\mathbb{Z})$. It is equipped with the structure of a projective variety.  

The result is often quoted by Teichm\"uller theorists in the following weak form: For a surface $S$ of genus $\geq 2$, if $\phi:S\to S$ is a holomorphic and bijective map that induces the identity on $H_1(S,\mathbb{Z}/n\mathbb{Z})$ for some $n\geq 3$, then $\phi$ is the identity map.

Grothendieck attributes this theorem to Serre. The result was announced in the first lecture (\cite{G1} , \S 2.4).

We now review the content of \cite{G10}.

  Grothendieck starts by saying that he will give the construction of Teichm\"uller space only in the case of genus $g\geq 2$. The construction in the case $g=1$ is left to the reader. This case is direct and elementary and does not use Hilbert schemes. He adds that this construction works more generally for moduli of tori and of polarized abelian varieties, but that it is transcendental. The construction that he gives for genus $g\geq 2$ is algebraic and leads to the existence of the moduli spaces (as schemes) of curves of genus $g$.

Grothendieck uses cohomology techniques which he did not use in the previous talks. He introduces the notion of \emph{linear rigidification}.\index{rigidification!linear}

Let $S$ be as before an analytic space, and $X/S$ a fiber space over $S$, where the fibers are Riemann surfaces. Let  $\Omega^1_{X/S}$ be  the sheaf  over $X$ whose sections are the holomorphic 1-forms on the fibers. These are called \emph{relative} 1-forms. The reason is that they only act on vectors that are tangent to the fibers.
  
We fix an integer $g\geq 2$ and an integer $k\geq 3$. For such a choice of $k$ and for any Riemann surface of genus $g$ over $S$, the tensor power $(\Omega^1_{X/S})^{\otimes k}$ is a relatively ample sheaf, that is, it has enough sections to allow embeddings of each fiber in a projective space. The dimension of the space of global 1-forms on the fibers is independent of the conformal structure of the surface. It can be computed using the Riemann-Roch theorem, and it depends only on the genus. 
 For every point on the Riemann surface, the space of all sections that vanish at this point is a codimension-one subspace of the vector space of cubic differentials.  Grothendieck uses this to define an embedding of the curve in the projective space. The image has high codimension.

\begin{definition}[Linear rigidification]\index{linear rigidification}\index{rigidification!linear}
Let $X$ be a curve of genus $g$ over an analytic space $S$. A \emph{linear rigidification}\index{linear rigidification} of $X$ is a sequence of sections $(w_1,\ldots,w_r)$ of the module $\mathcal{E}_{X/S}$ that are pointwise the base of the $k$-forms.\footnote{$r$ is determined by Rienamm-Roch $r=k(2g-2)-g+1$}

\end{definition}
Such a basis plays the role of a marking.\index{marking!sections}\index{marking!linear rigidification} This explains the word ``rigidification."\index{rigidification}

 For $s$ in $S$, the embedding of the fiber is obtained by mapping each element $p\in f^{-1}(s)$ to \[
 w_1(s)(p) :  w_2(s)(p) : \ldots  : w_r(s)(p) \in \mathbb{P}^{r-1}(\mathbb{C}).
\]
For $r=3$, this is the so-called tricanonical embedding of the fiber. It commutes with the symmetries of the curve: any symmetry of the curve is realized as a global projective motion of the projective space.

The set $ \mathcal{P}(\mathcal{E}_{X/S})$ of linear rigidifications of $X/S$ is an analytic variety over $S$. It is identified with the set of sections of a  locally trivial analytic principal fiber bundle whose structural group is $G=\mathrm{GL}(r,\mathbb{C})$. The automorphisms of the Riemann surface give rise to automorphisms of this group. Once we fix a basis for the vector space of cubic differentials, there are no more symmetries allowed in the group $G=\mathrm{GL}(3,\mathbb{C})$.

In conclusion, we get a canonical immersion \[X\to \mathcal{P}(\mathcal{E}_{X/S})\] which is functorial in $X/S$ (for a fixed $S$). 
This immersion is compatible with base-change and any automorphism of $X$ which induces the identity on $(\mathcal{E}_{X/S})$ is the identity of $X$. In other words, the morphism 
\[X\to \mathcal{P}(\mathcal{E}_{X/S})\] is a functor. It assigns to $X$ all the different ways of putting above $S$ a Riemann surface with a basis for the cubic $(k$-)differentials, everything varying holomorphically.

The set $\mathcal{R}(X/S)$ of linear rigidifications of $S$ is seen now as a $\mathrm{GL}(r,\mathbb{C})$-principal bundle. This is a locally trivial analytic principal fiber bundle over $S$.

The set of isomorphism classes of curves of genus $g$ above $S$ equipped with a linear rigidification is denoted by $\mathcal{U}_{\mathcal{R}}(S)$. The morphism 
\[
 \mathcal{U}_{\mathcal{R}} : S\to  \mathcal{U}_{\mathcal{R}}(S)
\]
is a contravariant functor 
from the category  $\mathbf{(An)}$ of analytic spaces to the category $\mathbf{(Ens)}$.

Now let $\mathcal{B}$ be the contravariant functor from
($\mathbf{An}$) into $\mathbf{(Ens)}$ which assigns to every Riemann surface $S$ the set $\mathcal{B}(S)$ of analytic closed subspaces of the projective space which are curves of genus $g$ above $X$. (This is the Hilbert functor,\index{functor!Hilbert} represented by the Hilbert scheme.) It was proved in Expos\'e IX (\cite{G9}, 1.3) that this functor is representable and this is where the Jacobi functor is used.
One needs to prove that the functor $\mathcal{U}_{\mathcal{R}}$ is representable relatively to $\mathcal{B}$ (\cite{G3}, IV, 3.3).

 Grothendieck proves the following result which is an analogue of the result in \cite{G1} that concerns the Teichm\"uller curve:
\begin{theorem}
The contravariant functor $\mathcal{U}_{\mathcal{R}}$ on the category ($\mathrm{A}_n$) is representable by an analytic space $M_{\mathcal{R}}$ all of whose local rings are regular of dimension $3g-3 + \mathrm{dim}(G)$. 
\end{theorem}

Using again the Riemann-Roch theorem, Grothendieck obtains the following 
\begin{corollary} [Corollary 3.3 of \cite{G10}]  Let $X$ be a curve of genus $g\geq 2$ over $S$. Then any $S$-automorphism of $X$ that induces the identity on the fibres of $X$ is the identity.
\end{corollary}
\begin{lemma} [Lemma 3.4 of \cite{G10}]  Let $X$ be an algebraic curve of genus $g\geq 2$ over a field $k$, let $A$ be its homogenous Jacobian variety and let $u$ be an automorphism of $X$ that induces the identity on $A$. Then $u$ is the identity.
\end{lemma}

\begin{lemma}[Lemma 3.5 of \cite{G10}]  Let $A$ be an abelian variety over an algebraically closed  field, let $n$ be an integer $\geq 3$ which is relatively prime to the characteristic of $k$, and let $u$ be a finite order automorphism of $A$ that fixes the points of order $n$ of $A$. Then $u$ is the identity.
\end{lemma}

The paper concludes with the following result on the universal Teichm\"uller curve:

\begin{corollary} [Corollary 5.4 of \cite{G10}]  The automorphism group of the universal Teichm\"uller curve of genus $g$ is reduced to the identity if $g\geq 3$ and it is the group $\mathbb{Z}/2\mathbb{Z}$ generated by the hyperelliptic symmetry in the cases $g=1$ and $g=2$.
\end{corollary}

Here, the universal automorphism group is the group of automorphisms of the Teichm\"uller universal curve which act  trivially on Teichm\"uller space.

Grothendieck concludes from the series of talks that we are reviewing that the analytic space $\mathcal{T}$ equipped with the automorphism group $\gamma$ can be considered as a satisfactory solution of the ``moduli problem" for curves of genus $g$, since it allows to reconstitute completely the fibered category $\mathcal{F}$ that we started with, and thus, to solve, at least theoretically, all the problems that can be expressed in terms of this fibered category. This statement is proved in \cite{G10}.

\section{In way of a conclusion} 

Our aim in this chapter was to present a point of view on Teichm\"uller theory which is completely ignored by Teichm\"uller theorists. In discussing with algebraic geometers and in corresponding with them, we also learned that this part of Grothendieck's work is also very poorly known in this community.

It is also a fact that algebraic geometers interested in moduli spaces have a very poor knowledge of most of the developements that were worked out by the analysts and the low-dimensional topologists on Teichm:"uller space.

 Let us also note that Grothendieck, in these lectures, raised several questions, among them the question of whether Teichm\"uller space is Stein (\cite{G1} p. 14).  In 1964, Bers and Ehrenpreis  showed that any Teichm\"uller space can be embedded as  a domain of holomorphy in some $\mathbb{C}^N$. 
  
  Finally, Grothendieck notes that it would be good to find a direct algebro-geometric description, in terms of universal problems, of the Baily-Satake compactifications of various moduli spaces (for curves as well as for polarized abelian varieties). He says that such a description would be applicable as well in the context of schemes, and would be the starting point of a purely geometric theory of automorphic functions, such as the one developed in \cite{Igusa} by Igusa for genus one. We shall not attempt to survey the developments of this rich theory, because this would lead us too far.

In conclusion, Teichm\"uller's 1944 paper and the corresponding Cartan seminar talks by Grothendieck that are considered in this  chapter shed a completely new light on Teichm\"uller space, and they present an essential aspect about families
of Riemann surfaces which is lacking in the Ahlfors-Bers point of view on this space, which is the most widely known.

\printindex

 \end{document}